\documentclass[pre,aps,twocolumn,showpacs]{revtex4}

\usepackage{epsf}
\usepackage{amsmath}
\usepackage{amsfonts}
\usepackage{amssymb}
\usepackage{bm}
\usepackage{bbm}
\usepackage{nicefrac}
\usepackage{color}
\usepackage{pifont}
\usepackage{graphicx}
\usepackage{epsfig}
\usepackage{latexsym}
\usepackage{pifont}

\DeclareMathOperator{\real}{\textrm{Re}}
\DeclareMathOperator{\imag}{\textrm{Im}}

\newcommand{\zi}{z_{\textrm{I}}}
\newcommand{\zr}{z_{\textrm{R}}}

\newcommand{\nui}{\nu_{\textrm{I}}}
\newcommand{\nur}{\nu_{\textrm{R}}}

\DeclareMathOperator{\erfc}{erfc}

\def\dd{{\mathrm{d}}}
\def\ii{{\mathrm{i}}}
\def\ee{{\mathrm{e}}}

\def\Ai{{\mathrm{Ai}}}

\def\tfrac#1#2{  {\textstyle{#1  \over #2}}   }

\def\Ai{   {\mathrm{Ai}}   }
\def\Bi{   {\mathrm{Bi}}   }

\renewcommand{\Re}{\mathrm{Re}}
\renewcommand{\Im}{\mathrm{Im}}

\begin{document}

\title{Numerical calculation of Bessel, Hankel and Airy functions}

\author{U. D. Jentschura}

\affiliation{Department of Physics,
Missouri University of Science and Technology,
Rolla, Missouri 65409-0640, USA}

\author{E. L\"{o}tstedt}

\affiliation{Department of Chemistry, School of Science, 
The University of Tokyo, 7-3-1 Hongo, Bunkyo-ku, Tokyo 113-0033, Japan}

\begin{abstract} 
The numerical evaluation of an individual Bessel or Hankel
function of large order and large argument is a notoriously problematic issue
in physics.  Recurrence relations are inefficient when an individual function of
high order and argument is to be evaluated. The coefficients in the well-known
uniform asymptotic expansions have a complex mathematical structure which
involves Airy functions. For Bessel and Hankel
functions, we present an adapted algorithm which relies on a 
combination of three methods: (i)
numerical evaluation of Debye polynomials, (ii) calculation of Airy functions
with special emphasis on their Stokes lines, and (iii) resummation of the
entire uniform asymptotic expansion of the Bessel and Hankel functions by
nonlinear sequence transformations.

In general, for an evaluation of a special function, we advocate the use of
nonlinear sequence transformations in order to bridge the gap between the
asymptotic expansion for large argument and the Taylor expansion for small
argument (``principle of asymptotic overlap'').  This general principle needs
to be strongly adapted to the current case, taking into account the complex
phase of the argument.  Combining the indicated techniques, we observe that it
possible to extend the range of applicability of existing algorithms. Numerical
examples and reference values are given.
\end{abstract}


\maketitle


\setlength{\parskip}{2.0ex plus0.8ex minus0.5ex}

%
%
\section{INTRODUCTION}

Bessel, Hankel and Airy functions constitute some of the most important special
functions used in theoretical physics, and their calculation is notoriously
problematic for extreme ranges of argument and order, even if their
mathematical definition is straightforward.  Especially, it should be noted
that recurrence relations are useful for arrays of Bessel or Hankel functions
when, for given argument, all functions up to a maximum order are needed.
However, recurrence relations cannot be used to good effect if an individual
function of large argument and order needs to be evaluated.

Notably, Bessel, Hankel and Airy functions occur in the multipole
decompositions of various operators in electrodynamics; these are known to be
slowly convergent decompositions in many cases. A lot of work has been
invested into the development of asymptotic expansions which may be used in the
calculation of the special functions.  After the famous paper of
Debye~\cite{De1909}, which used a saddle point expansion, a first treatise of
the notoriously problematic case of equal order and argument appeared in
Ref.~\cite{Ai1916}.  Further historical papers, where the theory was refined,
can be found in
Refs.~\cite{Ni1910,Wa1918,La1931,La1932,La1949,Ch1949,Ch1950,Gr1989}.  The
standard textbook~\cite{Wa1922} contains a collection of very useful formulas.
The known asymptotic formulas for large order (at fixed argument) are given in
reference volumes (e.g., Refs.~\cite{AbSt1972,OlLoBoCl2010,NISTlib}),
and they can be used,
together with asymptotic formulas for other Green functions~\cite{MoIn1995},
for the calculation of the properties of bound electrons.
Indeed, Bessel, Hankel and Airy functions belong to the most 
important special functions used in theoretical physics;
the much revived interest in these is also manifest in a recent 
monograph on Airy functions~\cite{VaSo2010}.
The renewed interest in the theory of special functions is 
also manifest in a number of other recent books and review 
articles~\cite{CuBrVeWaJo2008,GiSeTe2007,Te2007,GiSeTe2011}.

In a marvelous {\em tour de force}, Olver~\cite{Ol1954asymp1,Ol1954asymp2} has
derived uniform asymptotic expansions which hold for large order of Bessel and
Hankel functions, uniformly in the complex plane of the argument variable (for
arguments $z$ with a complex phase $|\arg(z) | < \pi - \epsilon$).  The
derivation is based on the general asymptotic properties of solutions of
second-order differential equations. The findings are summarized in Chapter~10
of the textbook~\cite{OlReprint} which is usually more accessible than the original
references (see also Chapter~8 of the recent Ref.~\cite{GiSeTe2007}).  

A central question surrounding the use of the uniform asymptotic expansions has
been their practical applicability to the calculation of Bessel and Hankel
functions. This question is important because the expansions, while uniformly
applicable in the complex plane, have a complicated mathematical structure, and
because they involve Airy functions whose numerical evaluation is eventually
required for arbitrary magnitude and complex phase of the argument. Despite
considerable and perhaps justified doubts regarding their usefulness for
numerical calculations, the uniform asymptotic
expansions~\cite{Ol1954asymp1,Ol1954asymp2} seem to be the most powerful
ones available for the calculation of an individual Bessel or Hankel function of
large argument and order. The aim of the current article is to show that their
domain of usefulness can be drastically enhanced if they are combined with the
``principle of asymptotic overlap'' that makes it possible to join asymptotic
regions for large argument with regions of small argument via the use of a
nonlinear sequence transformation in overlapping regions
[see also Section~2.4.1 of Ref.~\cite{CaEtAl2007}].

The importance of the numerical evaluation of an individual Bessel functions
for physics is highlighted by the substantial work devoted to the development
of asymptotic expansions and numerical
algorithms~\cite{Ba1966,ScAnGo1979,Te1994,Te1997uniform,HoOD1999,%
Pa2001hadamard1,Pa2001hadamard2,Pa2004hadamard3,ShWo2009,SmdHTi2009,BaMisc}.
If one aims to develop the algorithms based on the uniform asymptotic
expansions~\cite{Ol1954asymp1,Ol1954asymp2}, one first needs to evaluate Airy
$\Ai$ and $\Bi$ functions. For large modulus of the argument and variable
complex phase, their behavior is characterized by a Stokes phenomenon.  Along
the Stokes lines, i.e., along specific values of the complex phase of the
argument of the Airy functions, the relative magnitude the contribution of the
different saddle points changes.  A numerical algorithm for the Airy functions
has been described in~\cite{FaLoOl2004}.  It relies on a separation into the
region of large argument, where an asymptotic expansion is applied, and the
region of small argument, where a differential equation is being integrated.
Here, we advocate the use of a nonlinear sequence transformation in order to
bridge the gap between the asymptotic regime of large argument, and the regime
of small argument where a power series can be used.  The asymptotic
expansion used for large argument has to be be adapted according to the complex
phase of the argument.  As a byproduct of our analysis, we derive some higher
terms in the analytic expansion at the exact turning point $\nu = z$.

The paper is organized as follows: We first recall basic formulas in
Section~\ref{basic}.  Numerical calculations are described in
Section~\ref{summation}.  Analytic properties at the turning point $\nu = z$ are
calculated in Section~\ref{turning}.  Finally, conclusions are drawn in
Section~\ref{conclu}.
The Appendix~\ref{appa} is devoted to a general discussion about 
the saddle points in the complex plane related to the Bessel functions, and 
about the possibility of constructing an alternative algorithm.

%
%
\section{BASIC FORMULAS}
\label{basic}

For complex argument $z$ with $\Re(z) > 0$,
the evaluation of Bessel $J$ functions
can be traced to the evaluation of integrals
of the form 
[see Eq.~(10.9.6) of Ref.~\cite{OlLoBoCl2010}]
\begin{align} 
\label{defJ} 
& J_\nu(z) = \frac{1}{\pi} \int_0^\pi 
\cos\left[ z \, \sin(\theta) - 
\nu \, \theta \right] \, \dd \theta 
\nonumber\\[1ex]
& \quad - \frac{\sin(\nu \pi)}{\pi} \, 
\int_0^\infty 
\exp[- z \, \sinh(\theta) - \nu \, \theta] \, 
\dd \theta \,,
\end{align}
where the order $\nu$ of the Bessel function 
is not necessarily an integer
and $|\arg (z)| < \pi/2$.
Of particular interest are the Bessel $J$ functions,
as they are regular at the origin for positive integer $\nu$.
For integer $\nu$, the second term in the definition
of $J$ according to Eq.~\eqref{defJ} vanishes.
All of the definitions used here for Bessel functions,
and spherical Bessel functions, are contained in 
Chaps.~9 and~10 of Ref.~\cite{AbSt1972}.
Indeed, these and many of the asymptotic formulas 
used in the following are also included
in the modernized handbook~\cite{OlLoBoCl2010,NISTlib}.
Reference~\cite{AbSt1972} is now 
somewhat outdated but still the standard classic reference on 
the matter.

For half-integer $\nu$, the Bessel $J$ functions
are related to spherical Bessel functions according to the 
formula ($\ell$ is an integer)
\begin{equation}
j_\ell(z) = \sqrt{\frac{\pi}{2 \,z}} \, J_{\ell + 1/2}(z) \,.
\end{equation}
This relation is given in Eq.~(10.47.3) of Ref.~\cite{OlLoBoCl2010}.
For $\Re(z) > 0$, the Bessel $Y$ function is defined as
[see Eq.~(10.9.7) of Ref.~\cite{OlLoBoCl2010}]
\begin{align} 
\label{defY} 
& Y_\nu(z) = \frac{1}{\pi} 
\int_0^\pi \sin\left[ z \, \sin(\theta) - 
\nu \, \theta \right] \, \dd \theta 
\nonumber\\[1ex]
& \quad - \frac{1}{\pi} \, 
\int_0^\infty 
\left\{ \ee^{\nu \,t} + 
\ee^{-\nu \, t} \cos(\nu \, \pi) \right\} \, 
\ee^{ - z \, \sinh(t) } \, \dd \theta \,,
\end{align}
for $|\arg (z)| < \pi/2$.
The spherical $y$ function is defined as
\begin{equation}
y_\ell(z) = \sqrt{\frac{\pi}{2 \,z}} \, Y_{\ell + 1/2}(z) \,.
\end{equation}
This relation is given in Eq.~(10.47.4) of Ref.~\cite{OlLoBoCl2010}.
The Hankel functions are defined as
[see Eqs.~(9.1.3) and~(9.1.4) of Ref.~\cite{AbSt1972}]
\begin{subequations}
\begin{align}
\label{defH1}
H^{(1)}_\nu(z) =& \; J_\nu(z) + \ii \, Y_\nu(z) \,,
\\[1ex]
\label{defH2}
H^{(2)}_\nu(z) =& \; J_\nu(z) - \ii \, Y_\nu(z) \,,
\\[1ex]
\label{defh1}
h^{(1)}_\ell(z) =& \; j_\ell(z) + \ii \, y_\ell(z) \,,
\\[1ex]
\label{defh2}
h^{(1)}_\ell(z) =& \; j_\ell(z) - \ii \, y_\ell(z) \,.
\end{align}
\end{subequations}
The definitions of the spherical Hankel functions are given 
in Chap.~10.1.1 of Ref.~\cite{OlLoBoCl2010}.

The integral representations~\eqref{defJ} and~\eqref{defY}
are valid for $\Re (z) > 0$. For purely imaginary $z$, we may use 
a definition in terms of the modified Bessel functions 
$I_\nu (x)$ and $K_\nu (x)$, see 
Eqs.~\eqref{eq:modifiedIK1}--\eqref{eq:modifiedIK4}. 
Below, we describe a numerical algorithm
with the notion $0 \leq \arg(z) < \pi$
in mind. Arguments $z$ with $-\pi \leq \arg(z) < 0$
are treated by the transformation $z \to z' = z \, \exp(\ii \, \pi)$
(so that the transformation $z \to z'$ does not leave 
the first Riemann sheet).
The conversion formulas can be derived based on 
Eqs.~(9.1.35)---(9.1.39) of Ref.~\cite{AbSt1972} and read
\begin{subequations}
\label{six}
\begin{align}
J_\nu(z) =& \; 
\ee^{-\ii \, \nu \, \pi} \, J_\nu(z \, \ee^{\ii \, \pi} ) \,,
\quad
-\pi \leq \arg(z) < 0 \,,
\\[1ex]
Y_\nu(z) =& \; \ee^{\ii \, \nu \, \pi} \, 
Y_\nu(z \, \ee^{\ii \, \pi} ) -
2 \, \ii \, \cos(\nu \, \pi) \, J_\nu(z \, \ee^{\ii \, \pi} ) \,,
\\[1ex]
H^{(1)}_\nu(z) =& \; 
\ee^{\ii \, \nu \, \pi} \, H^{(1)}_\nu(z \, \ee^{\ii \, \pi}) + 
2\, \ee^{-\ii \, \nu \, \pi} \, J_\nu(z \ee^{\ii\pi}) \,,
\\[1ex]
H^{(2)}_\nu\left( z \right) =& \; 
- \ee^{\ii \, \nu \, \pi} \, H^{(1)}_\nu(z \, \ee^{\ii \, \pi} ) \,.
\end{align}
\end{subequations}
(Many of the asymptotic formulas given here are also included
in the modernized handbook~\cite{OlLoBoCl2010}, but for the time being,
we prefer to refer to equation references in the somewhat outdated,
but standard classic Ref.~\cite{AbSt1972}.)
Alternatively, one may use direct complex conjugation,
\begin{subequations}
\label{seven}
\begin{align}
J_\nu(z) =& \; 
\left( J_\nu(z^*) \right)^* \,,
\quad
Y_\nu(z) =
\left( Y_\nu(z^*) \right)^* \,,
\\[1ex]
H^{(1)}_\nu(z) =& \; 
\left( H^{(2)}_\nu(z^*) \right)^* \,,
\;
H^{(2)}_\nu(z) =
\left( H^{(1)}_\nu(z^*) \right)^* \,,
\end{align}
\end{subequations}
where $z^*$ is the complex conjugate of $z$.
Using a combination of the formulas~\eqref{six} and~\eqref{seven},
we could in principle restrict the range of complex phases
of the arguments to the first quadrant of the complex $z$ plane.
However, as evident from Eqs.~\eqref{Juni},~\eqref{Yuni},
\eqref{H1uni} and~\eqref{H2uni} below,
we would still need to evaluate the Airy $\Ai$ and $\Bi$ 
functions in the entire complex plane, even if we
restrict $z$ to the first quadrant in the 
complex $z$ plane (and the former constitutes the 
main computational challenge). A simple restriction to the 
upper half of the complex $z$ plane thus seems to be 
most effective.

Without loss of generality, we restrict our attention to the case 
$\nu > 0$ in the following. For $\nu < 0$, the conversion 
formulas are as follows,
\begin{subequations}
\label{conv}
\begin{align}
\label{conva}
H^{(1)}_{-\nu}(z) =& \; \ee^{\ii \, \pi \nu} \, H^{(1)}_{\nu}(z)  \,,
\\[1ex]
\label{convb}
H^{(2)}_{-\nu}(z) =& \; \ee^{-\ii \, \pi \nu} \, H^{(2)}_{\nu}(z) \,,
\\[1ex]
\label{convc}
J_{-\nu}(z) =& \; \cos(\pi \nu) \, J_{\nu}(z) - \sin(\pi \nu) \, Y_{\nu}(z) \,,
\\[1ex]
\label{convd}
Y_{-\nu}(z) =& \; \sin(\pi \nu) \, J_{\nu}(z) + \cos(\pi \nu) \, Y_{\nu}(z) \,.
\end{align}
\end{subequations}
The conversion matrix for the Bessel functions
$J$ and $Y$ has the same structure as a rotation matrix
for an angle $\pi \nu$. 
For the Hankel functions, the above formulas~\eqref{conva}
and~\eqref{convb} can be found in
Eqs.~(9.1.5) and (9.1.6) of Ref.~\cite{AbSt1972}.

In principle, one might speculate that 
the above integral representations~\eqref{defJ}
and~\eqref{defY} should be
sufficient in order to 
numerically evaluate an individual Bessel function.
However, the numerical difficulties for large $\nu$ are 
nearly insurmountable in view of apparent numerical
oscillations of the integrand.
While one can investigate complex integration contours 
with the notion of adopting a 
steepest descent method
(see Appendix~\ref{appa}),
these representations do not immediately lead
to a uniformly applicable algorithm, either.

Finally, let us recall the basic 
asymptotic properties of Bessel $J$ and $Y$ functions
for $\nu > 0$. Only the Bessel $J$ functions is
regular at the origin, and we have
\begin{subequations}
\begin{align}
J_\nu(z) \sim & \; 
\frac{ 1 }{ \Gamma(\nu + 1) } \,
\left( \frac{z}{2} \right)^\nu \,,
\quad |z| \to 0 \,,
\\[1ex]
Y_\nu(z) \sim & \; -\frac{ \Gamma(\nu) }{ \pi } \,
\left( \frac{2}{z} \right)^\nu \,, \quad |z| \to 0 \,.
\end{align}
\end{subequations}
The literature on Bessel functions is manifold.
A very useful reference is the standard treatise~\cite{Wa1922}.
In Chapter~10 of Ref.~\cite{OlReprint}, basic asymptotic expansions and
properties of Bessel $J$ and $Y$ functions, and of their
derivatives, are reviewed and explained very clearly.

One is often faced with the problem of calculating 
strings of Bessel functions whose indices 
differ by integers~\cite{Mo1974a,Mo1974b,LoJe2009}.
Recursive algorithms based on the relations
\begin{subequations}
\label{recur}
\begin{align}
J_{\nu - 1}(x) + J_{\nu + 1}(x) = &\;
\frac{2 \nu}{x} \, J_\nu(x) \,,
\\[1ex]
Y_{\nu - 1}(x) + Y_{\nu + 1}(x) = &\;
\frac{2 \nu}{x} \, Y_\nu(x) \,,
\end{align}
\end{subequations}
can be very effective, as explained in Section~10.5 on p.~452 of
Ref.~\cite{AbSt1972}.  For Bessel $J$ functions, one starts a three-term
downward recursion in $\nu$ with two essentially arbitrary starting values for
$J_{\nu + 1}(x)$ and $J_\nu(x)$ at high $\nu$ and continues to calculate
$J_{\nu-1}(x)$ until the order of the Bessel 
function becomes zero.  One can then either calculate
$J_0(x)$ explicitly and use the recurrence relation
upwards (filling the array of Bessel functions), 
or fix the normalization of all calculated Bessel
functions by a normalization
condition~\cite{Mi1950,BiEtAl1960,Ga1967}
(see also Chap.~10.10.5 of Ref.~\cite{AbSt1972}). In Ref.~\cite{Ga1967}, the
computational aspects of three-term recursion relations have been discussed
with a special emphasis on their numerical stability.  Here, we are dealing
with a different problem, namely, the evaluation of an individual Bessel
function of high order and argument $J_\nu(x)$ and $Y_\nu(x)$, without recourse
to any recurrence relation in $\nu$.

%
%
\section{SUMMATION OF THE UNIFORM ASYMPTOTICS}
\label{summation}

%
%
\subsection{Uniform asymptotic expansions}

The task in the current investigation 
is to calculate the functions
\begin{subequations}
\begin{align}
J_\nu(z = \nu \, y) \,, \qquad
J'_\nu(z = \nu \, y) \,, 
\\[1ex]
Y_\nu(z = \nu \, y) \,, \qquad
Y'_\nu(z = \nu \, y) 
\end{align}
\end{subequations}
for complex argument $-\pi \leq \arg(z) < \pi$, and real $\nu > 0$. 
In a numerical code, it is sufficient to treat the 
complex phase range $0 \leq \arg(z) < \pi$.
The range $-\pi < \arg(z) < 0$ is covered by 
Eqs.~\eqref{six} and~\eqref{seven}.
Because we are using asymptotic expansions valid for
large $\nu$, we also assume that $\nu > 50$.
For $\nu < 50$, one may use Miller's 
method~\cite{Mi1950,BiEtAl1960,Ga1967}.
A brief digression on this algorithm can also 
be found in Chap.~10.10.5 of Ref.~\cite{AbSt1972}.
The case $\nu < -50$ then is covered by Eq.~\eqref{conv}.
The parameterization $z = \nu \, y$ 
is useful to identify the notoriously problematic
region near $y \approx 1$.
It is also being used below in Appendix~\ref{appa}.

We first have to recall the
uniform asymptotics of the Bessel
functions (see Refs.~\cite{Ol1954asymp1,Ol1954asymp2}),
These are also listed in Eqs.~(9.3.35),~(9.3.36),~(9.3.43) 
and~(9.3.44) of Ref.~\cite{AbSt1972},
and in Chapter~10 of Ref.~\cite{OlReprint}.
A brief rederivation is 
given in Ref.~\cite{HoOD1999}.
The uniform asymptotics for the Bessel $J$ function are given by
\begin{align}
\label{Juni}
J_\nu(\nu \, y) \sim & \;
\left( \frac{4 \zeta}{1 - y^2} \right)^{1/4}  \, \left\{ 
\frac{ \Ai(\nu^{2/3} \, \zeta) }{\nu^{1/3}} \;
\sum_{k=0}^\infty \frac{ a_k(\zeta) }{ \nu^{2 k} }  \right.
\nonumber\\[1ex]
& \left. + \frac{ \Ai'(\nu^{2/3} \, \zeta) }{\nu^{5/3}} \;
\sum_{k=0}^\infty \frac{ b_k(\zeta) }{ \nu^{2 k} } \right\} \,.
\end{align}
This asymptotic formula is valid for 
$\nu \to \pm \infty$ and
$\arg (y) \leq \pi - \epsilon$,
where $\epsilon$ is an arbitrarily small positive number.
We denote the Airy function of the first kind as~$\Ai$.

For $y \geq 0$, the $\zeta$ variable is defined as
\begin{subequations}
\begin{align}
\label{a}
\zeta =& \; \left( \frac{3}{2} \right)^{2/3} \,
\left[ \ln\left( \frac{1 + \sqrt{1 - y^2}}{y} \right) - 
\sqrt{1 - y^2} \right]^{2/3} > 0 \,, 
\nonumber\\[1ex]
& \; 0 \leq y < 1 \,,
\\[1ex]
\label{b}
\zeta =& \; - \left( \frac{3}{2} \right)^{2/3} \,
\left[ \sqrt{y^2 - 1}  - \arccos\left( \frac{1}{y} \right) \right]^{2/3} < 0 \,,
\nonumber\\[1ex]
& \; y > 1 \,.
\end{align}
\end{subequations}
The calculation of $\zeta$ for complex $y$ relies on the 
formula
\begin{equation}
\frac23 \, \zeta^{3/2} = 
\ln\left( \frac{1 + \sqrt{1 - y^2}}{y} \right) - \sqrt{1 - y^2}
\end{equation}
where the branches take their principal values when 
$z \in (0,1)$ and $\zeta \in (0,\infty)$ and $\zeta$ is 
continuous elsewhere, as described in Chapter~10.1
of Ref.~\cite{OlReprint}.

The $a_k$ and $b_k$ coefficients entering 
Eqs.~\eqref{Juni}---\eqref{H2uni} read
\begin{subequations}
\begin{align}
\label{ak}
a_k(\zeta) = & \; \sum_{s = 0}^{2 k} \mu_s \; \zeta^{-3 s/2} \;
u_{2 k - s}[ (1 - y^2)^{-1/2} ] \,,
\\[1ex]
\label{bk}
b_k(\zeta) = & \; - \zeta^{-1/2} \;
\sum_{s = 0}^{2 k + 1} \lambda_s \; \zeta^{-3 s/2} \;
u_{2 k - s + 1}[ (1 - y^2)^{-1/2} ] \,.
\end{align}
\end{subequations}
They involve the Debye $u$ polynomials,
and coefficients $\mu_s$ and $\lambda_s$
which need to be defined.
The corresponding formulas read
\begin{subequations}
\begin{align}
\label{lambdas}
\lambda_s =& \; \frac{1}{s! \, 144^s}
\prod_{
\begin{array}{c} 
\scriptstyle 
m=2 s + 1 \\[-0.3ex]
\scriptstyle 
m \; \rm{odd} \end{array}}^{6 s - 1 } 
\!\!\!\! (m) 
\;\; = \;\; \frac{\Gamma(3 s + \tfrac{1}{2})}{9^s \, \sqrt{\pi} \, 
\Gamma\!\left(2 s + 1\right)} 
\nonumber\\
= & \;
\frac{1}{s! \, 144^s} 
(2s+1) \, (2s+3) \, \cdots \, (6s-1) \,,
\\[1ex]
\label{mus}
\mu_s =& \; - \frac{6 s + 1}{6 s - 1} \, \lambda_s =
- \frac{2 \, \Gamma(3 s + \tfrac{3}{2})}{9^s \, \sqrt{\pi} \, 
(6 s - 1) \, \Gamma\!\left(2 s + 1\right)} \,.
\end{align}
\end{subequations}
The Debye polynomials fulfill
$u_0(t) = 1$ and are otherwise defined recursively as 
\begin{equation}
u_{k+1}(t) = \tfrac12 \, t^2 \, (1-t^2) \, u'_k(t) 
+ \frac18 \, \int_0^t \dd t' \, (1 - 5 t'^2) \, u_k(t') \,.
\end{equation}
For polynomials, the operations 
of differentiation and integration
can be represented by simple multiplication operations
acting on a coefficient matrix.
This is due to the trivial identity $\dd x^n/\dd x =
n \, x^{n-1}$, applied to integer $n$.
On a computer system, it is thus possible to evaluate
the coefficients of, say, the polynomial coefficients 
for the first few hundred Debye
polynomials and to use them in order to evaluate
the $a_k(\zeta)$ and $b_k(\zeta)$ coefficients
for given $\zeta$.

The asymptotic expansion~\eqref{Juni}
obviously is an expansion for large $\nu$,
and it is valid even in the problematic region $y \approx 1$.
We are now in the position to give
the corresponding formula for the $Y$ function,
which involves the Airy function of the second kind $\Bi$ 
and its derivative,
\begin{align}
\label{Yuni}
Y_\nu(\nu \, y) \sim & \;
-\left( \frac{4 \zeta}{1 - y^2} \right)^{1/4}  \, \left\{ 
\frac{ \Bi(\nu^{2/3} \, \zeta) }{\nu^{1/3}} \;
\sum_{k=0}^\infty \frac{ a_k(\zeta) }{ \nu^{2 k} } 
\right.
\nonumber\\[1ex]
& \; \left. +
\frac{ \Bi'(\nu^{2/3} \, \zeta) }{\nu^{5/3}} \;
\sum_{k=0}^\infty \frac{ b_k(\zeta) }{ \nu^{2 k} } 
\right\}\,.
\end{align}
The uniform asymptotic expansion of the
Hankel $H^{(1)}$ function is given by
\begin{align}
\label{H1uni}
& H^{(1)}_\nu(\nu \, y) \sim 
2 \, \ee^{-\pi \ii /3} \,
\left( \frac{4 \zeta}{1 - y^2} \right)^{1/4}  \, 
\nonumber\\[1ex]
& \; \times \left\{ 
\frac{ \Ai(\ee^{2 \pi \ii/3} \, \nu^{2/3} \, \zeta) }{\nu^{1/3}} \;
\sum_{k=0}^\infty \frac{ a_k(\zeta) }{ \nu^{2 k} } 
\right. 
\nonumber\\[1ex]
& \; \qquad \left. +
\frac{ \ee^{2 \pi \ii/3} \; \Ai'(\ee^{2 \pi \ii/3} \, 
\nu^{2/3} \, \zeta) }{\nu^{5/3}} \;
\sum_{k=0}^\infty \frac{ b_k(\zeta) }{ \nu^{2 k} } 
\right\} \,.
\end{align}
According to Eq.~\eqref{defH1}, the 
uniform asymptotic expansion for $H^{(2)}$ is obtained by 
changing the sign of the imaginary unit,
\begin{align}
\label{H2uni}
& H^{(2)}_\nu(\nu \, y) \sim
2 \, \ee^{\pi \ii /3} \,
\left( \frac{4 \zeta}{1 - y^2} \right)^{1/4}  \, 
\nonumber\\[1ex]
& \; \times \left\{ 
\frac{ \Ai(\ee^{-2 \pi \ii/3} \, \nu^{2/3} \, \zeta) }{\nu^{1/3}} \;
\sum_{k=0}^\infty \frac{ a_k(\zeta) }{ \nu^{2 k} }  \right.
\nonumber\\[1ex]
& \; \qquad \left. +
\frac{ \ee^{-2 \pi \ii/3} \; \Ai'(\ee^{-2 \pi \ii/3} \, 
\nu^{2/3} \, \zeta) }{\nu^{5/3}} \;
\sum_{k=0}^\infty \frac{ b_k(\zeta) }{ \nu^{2 k} } 
\right\} \,.
\end{align}
The corresponding equations for the 
derivatives of the 
Bessel and Hankel functions can be found in Eqs.~(9.3.43),~(9.3.44)
and (9.3.45) of Ref.~\cite{AbSt1972}.
Otherwise, the derivatives are also accessible via the 
formula
\begin{equation}
\frac{\partial}{\partial z}{\cal J}'_\nu(z) \equiv
{\cal J}'_\nu(z) =
\frac12 \, 
\left( {\cal J}_{\nu-1}(z) - {\cal J}_{\nu+1}(z) \right) \,,
\end{equation}
where ${\cal J}$ stands for $J$, $Y$, $H^{(1)}$ or $H^{(2)}$.

%
%
\begin{figure}[t!]
\begin{center}
\includegraphics[width=0.6\linewidth]{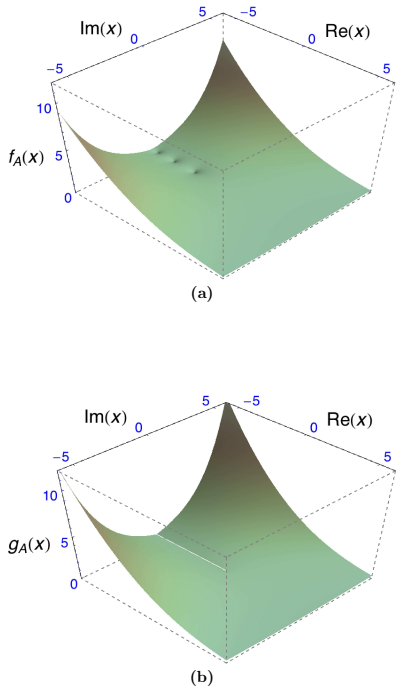}\\
\end{center}
\caption{\label{fig1} (Color online.) Figure (a) shows a contour plot of
$f_A(x) = |\Ai(x)|^{1/6}$ as a function of $\Re(x)$ and $\Im(x)$.
The exponent
$\tfrac16$ is introduced in order to prevent ``overflow'' of the 
plotted function near the boundaries of the considered range of arguments.
Figure~(b) shows a contour plot of $|\exp(- \tfrac23 \, x^{3/2})
|^{1/6}$ as a function of $\Re(x)$ and $\Im(x)$. For $x < 0$, i.e.~on the
negative real axis, the modulus is unity (Stokes line).  The zeros of $\Ai(x)$
give rise to the visible ``bump holes'' on the negative real axis in panel~(a).
Except for the region near $\arg(x) = \pm \pi$, the Airy $\Ai$ function can be
described using a single, uniform asymptotic formula $A(x)$ as defined in
Eq.~\eqref{defAx}, which is proportional to $|\exp(- \tfrac23 \, x^{3/2})
|^{1/6}$.}
\end{figure}

%
%
\begin{figure}[t!]
\begin{center}
\includegraphics[width=0.6\linewidth]{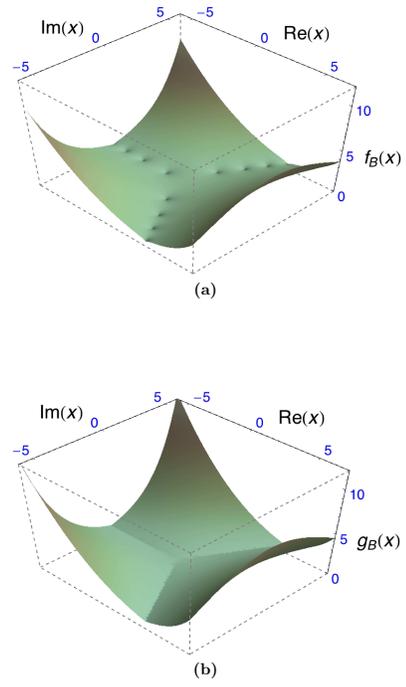}\\
\end{center}
\caption{\label{fig2} (Color online.) Figure (a) shows a plot of
$|\Bi(x)|^{1/6}$ as a function of $\Re(x)$ and $\Im(x)$. From the contour plot,
it is evident that a simple exponential of the form $|\exp(\pm \tfrac23 \,
x^{3/2})|^{1/6}$ cannot possibly describe the asymptotic behavior of the Airy
$\Bi$ function.  The zeros of $\Bi(x)$ give rise to visible ``bump holes'' on the
negative real axis and along the lines $\arg(x) = \pm \pi/3$ in panel~(a).
Figure~(b) shows a contour plot of $|f(x)|$ where $f(x) = |\exp(\tfrac23 \,
x^{3/2}) |^{1/6}$ for $| \arg(x) | \leq \pi/3$ and $f(x) = |\exp(-\tfrac23 \,
x^{3/2}) |^{1/6}$ for $\pi/3 < | \arg(x) | \leq \pi$, representing the Stokes
line behavior near $\arg(x) = \pm \pi/3$ and $\arg(x) = \pi$ [see
Eq.~\eqref{BiAsymp}].}
\end{figure}

%
%
\subsection{Evaluation of the Airy function}
\label{airy}

We now discuss the principle of asymptotic overlap
in the evaluation of the Airy functions of the 
first and second kind.
In contrast to the usual notation, we denote the complex 
argument of the Airy functions as $x$, in order to 
distinguish it from the argument $z$ of the Bessel 
and Hankel functions.
For $x \to 0$, we use the expansion,
\begin{subequations}
\label{AiBiPower}
\begin{align}
\label{AiSmall}
\Ai(x) = & \;
\sum_{k = 0}^\infty \frac{3^{-2k-\tfrac23}}%
{k! \, \Gamma\left(k + \tfrac23\right)} \, x^{3 \, k} 
- \sum_{k = 0}^\infty \frac{3^{-2k-\tfrac43} }%
{k! \, \Gamma\left(k + \tfrac43\right)} \, x^{3 \, k + 1} \,.
\end{align}
For the derivative of the Airy function, we have
\begin{align}
\label{ApSmall}
\Ai'(x) = & \;
-\sum_{k = 0}^\infty \frac{3^{-2k-\tfrac13}}%
{k! \, \Gamma(k + \tfrac13)} \, x^{3 \, k} 
+ \sum_{k = 0}^\infty \frac{3^{-2 k - \tfrac53}}%
{k! \, \Gamma(k + \tfrac53)} \, x^{3 \, k + 2}  \,.
\end{align}
The Airy $\Bi$ function is given by
\begin{align}
\label{BiSmall}
\Bi(x) = & \;
\sum_{k = 0}^\infty \frac{ 3^{-2k-\tfrac16} }%
{k! \, \Gamma(k + \tfrac23)} \, x^{3 \, k} 
+ \sum_{k = 0}^\infty \frac{3^{-2k-\tfrac56} }%
{k! \, \Gamma(k + \tfrac43)} \, x^{3 \, k + 1}  \,,
\end{align}
and its derivative by
\begin{align}
\label{BpSmall}
\Bi'(x) = & \;
\sum_{k = 0}^\infty \frac{3^{-2 k + \tfrac16} }%
{k! \, \Gamma(k + \tfrac13)} \, x^{3 \, k} 
+ \sum_{k = 0}^\infty \frac{3^{-2k-\tfrac76} }%
{k! \, \Gamma(k + \tfrac53)} \, x^{3 \, k + 2} \,.
\end{align}
\end{subequations}
The convergence radius of the expansions~\eqref{AiBiPower} is actually 
infinite, because of the Gamma functions in the denominator.
Still, they are faced with numerical problems 
for large negative $x$. 
In order to illustrate this fact, we draw an analogy
to the (likewise convergent) expansion $\exp(-x) = \sum_{k=0}^\infty 
(-1)^k x^k/k!$,
which also involves a power in the numerator and a factorial
in the denominator. Being absolutely convergent,
this expansion is numerically useless for the evaluation of
$\exp(-x)$ at $x = 5000$, because the largest term 
in the series is of order $\exp(5000)$, 
whereas the entire sum amounts to 
$\exp(-5000) \approx 3.37 \times 10^{-2172}$.
If the series were summed term by term, one would incur a
numerical loss of more than 4000~decimals.
Other methods (asymptotic expansions for large argument) therefore
have to be pursued.
Figures~\ref{fig1} and~\ref{fig2}
show the Stokes phenomenon. From Fig.~\ref{fig2},
we infer that the Airy $\Bi$ integral
is exponentially growing for $|x| \to \infty$,
along the complex directions
$\arg(x) = 0$ and $\arg(x) = \pm 2 \pi/3$,
with Stokes lines at $\arg(x) = \pm \pi/3$ and 
$\arg(x) = \pi$.
Therefore, the $\Bi$ integral cannot be represented by a
simple asymptotic divergent series but must be the sum of two.
This can be justified on the basis of saddle point considerations
\cite{Ol1954asymp1}.

The expansions~\eqref{AiBiPower} are valid for $x \to 0$. 
Complementing these expansions,
for $x \to \infty$, one is interested in 
suitable asymptotic expansions of the Airy functions.
To this end, it is useful to relate the Airy $\Ai$
and $\Bi$ functions to a modified Bessel 
functions $K$ and $I$ of order $\pm \tfrac{1}{3}$,
\begin{subequations}
\label{AiBiKI}
\begin{align}
\label{AiK}
\Ai(x) =& \; \frac{1}{\pi} \, \sqrt{\frac{x}{3}} \, 
K_{\tfrac{1}{3}}\!\!\left( \frac23 \, x^{3/2} \right),
\\[1ex]
\label{BiI}
\Bi(x) =& \; \sqrt{\frac{x}{3}} \, 
\left[ I_{\tfrac{1}{3}}\left( \frac23 \, x^{3/2} \right) +
I_{-\tfrac{1}{3}}\left( \frac23 \, x^{3/2} \right)
\right] \,.
\end{align}
\end{subequations}
Based on the relationship of the Airy functions
to modified Bessel functions,
we infer that the following asymptotic series
are relevant for the calculation of the Bessel functions at
large argument,
\begin{subequations}
\label{defABCD}
\begin{align}
\label{defAx}
A(x) =& \;
\frac{\exp\left(- \tfrac{2}{3} \, x^{3/2} \right)}{ \pi^{3/2} \, x^{1/4}}
\nonumber\\
& \; \times \sum_{k = 0}^\infty 
\frac{ (-3)^k \, \Gamma\left( k + \tfrac{1}{6} \right) \, 
\Gamma\left( k + \tfrac{5}{6} \right) }%
{2^{2 k + 2} \, k! \;\; x^{3 k / 2} } \,, 
\\[1ex]
\label{defBx}
B(x) =& \;
\frac{\exp\left(\tfrac{2}{3} \, x^{3/2} \right)}{ \pi^{3/2} \, x^{1/4}} 
\nonumber\\
& \; \times \sum_{k = 0}^\infty 
\frac{ 3^k \, \Gamma\left( k + \tfrac{1}{6} \right) \, 
\Gamma\left( k + \tfrac{5}{6} \right) }%
{2^{2 k + 2} \, k! \;\; x^{3 k / 2} } \,, 
\\[1ex]
\label{defCx}
C(x) =& \;
\frac{x^{1/4} \, \exp\left(-\tfrac{2}{3} \, x^{3/2} \right)}{ \pi^{3/2} } \,
\nonumber\\
& \; \times \sum_{k = 0}^\infty 
\frac{ (-3)^k \, \Gamma\left( k - \tfrac{1}{6} \right) \, 
\Gamma\left( k + \tfrac{7}{6} \right) }%
{2^{2 k + 2} \, k! \;\; x^{3 k / 2} } \,, 
\end{align}
whereas $D(x)$ is the only function that carries an overall 
minus sign in the prefactor,
\begin{align}
\label{defDx}
D(x) =& \;
-\frac{x^{1/4} \, \exp\left(\tfrac{2}{3} \, x^{3/2} \right)}{ \pi^{3/2} } \,
\nonumber\\
& \; \times \sum_{k = 0}^\infty 
\frac{ 3^k \, \Gamma\left( k - \tfrac{1}{6} \right) \, 
\Gamma\left( k + \tfrac{7}{6} \right) }%
{2^{2 k + 2} \, k! \;\; x^{3 k / 2} } \,.
\end{align}
\end{subequations}
These series are generalized hypergeometric series 
of the ${}_2 F_0$ form that
diverge for every nonzero argument $x^{-3/2}$ unless they terminate.
Their asymptotic relation
to the Airy functions is clarified below.
Indeed, the asymptotic expansions of $\Ai$ and $\Bi$ and of their derivatives
can be related to the series $A(x)$, $B(x)$, $C(x)$ and $D(x)$.
For $|x| \to \infty$, we have the divergent asymptotic expansion
as a function of the complex phase,
\begin{subequations}
\label{AiAsymp}
\begin{align}
\Ai(x) \sim & \; A(x) \,, 
\nonumber\\
& \; |x| \to \infty \,, \quad 
|\arg(x)| \leq \pi - \epsilon\,,
\\[1ex]
\Ai(x) \sim & \; A(x) + \ii \, B(x) \,, 
\nonumber\\
& \; |x| \to \infty \,, \quad \pi - \epsilon < \arg(x) \leq \pi \,,
\\[1ex]
\Ai(x) \sim & \; A(x) - \ii \, B(x) \,, 
\nonumber\\
& \; |x| \to \infty \,, \quad -\pi < \arg(x) < -\pi +\epsilon\,,
\end{align}
\end{subequations}
For the derivative of the Ai function, we have
\begin{subequations}
\label{ApAsymp}
\begin{align}
\Ai'(x) \sim & \; C(x) \,, 
\nonumber\\
& \; |x| \to \infty \,, \quad | \arg(x) | \leq \pi - \epsilon \,,
\\[1ex]
\Ai'(x) \sim & \; C(x) + \ii \, D(x) 
\nonumber\\
& \; |x| \to \infty \,, \quad \pi - \epsilon < \arg(x) \leq \pi \,,
\\[1ex]
\Ai'(x) \sim & \; C(x) - \ii \, D(x) 
\nonumber\\
& \; |x| \to \infty \,, \quad -\pi < \arg(x) < -\pi +\epsilon\,,
\end{align}
\end{subequations}
for arbitrarily small $\epsilon$.
The role of the $\epsilon$ parameter in these expressions
can be clarified as follows. From the formulas, it is evident that 
$\arg(x) \to \pm \pi$ (from below or above, respectively),
there is an admixture of the asymptotic $B$ and $D$ series 
to the $A$ and $C$ series. The magnitude of that admixture
is related to the magnitude of $|x|$. For large $|x|$, 
we can choose $\epsilon$ to be small. For finite $\epsilon$, 
one has to compare the magnitude of the 
two contributions.
If the subdominant saddle point is 
numerically significant, one has to add its contributions
(see also Fig.~\ref{fig3} below).

The dominant contributions to the asymptotic
expansions for the Airy $\Bi$ functions are given by
\begin{subequations}
\label{BiAsymp}
\begin{align}
\Bi(x) \sim & \; 2 \, B(x) \,, 
\nonumber\\
& \; |x| \to \infty \,, \quad 
|\arg(x)| < \frac{\pi}{3} - \epsilon \,,
\\[1ex]
\Bi(x) \sim & \; 2 B(x) + \ii \, A(x) \,, 
\nonumber\\
& \; |x| \to \infty \,, \quad 
\frac{\pi}{3} - \epsilon < \arg(x) < \frac{\pi}{3} + \epsilon \,,
\\[1ex]
\Bi(x) \sim & \; \ii \, A(x) \,, 
\nonumber\\
& \; |x| \to \infty \,, \quad 
\frac{\pi}{3} + \epsilon < \arg(x) < \pi - \epsilon\,,
\\[1ex]
\Bi(x) \sim & \; B(x) + \ii \, A(x) \,, 
\nonumber\\
& \; |x| \to \infty \,, \quad 
\pi - \epsilon < \arg(x) \leq \pi \,.
\end{align}
For $\arg(x) < 0$, the corresponding conditions 
are obtained by complex conjugation,
\begin{align}
\Bi(x) \sim & \; 2 B(x) - \ii \, A(x) \,, 
\nonumber\\
& \; |x| \to \infty \,, \quad 
-\frac{\pi}{3} -\epsilon < \arg(x) < -\frac{\pi}{3} + \epsilon \,,
\\[1ex]
\Bi(x) \sim & \; -\ii \, A(x) \,, 
\nonumber\\
& \; |x| \to \infty \,, \quad 
-\frac{\pi}{3} + \epsilon < \arg(x) < -\pi + \epsilon \,,
\nonumber\\[1ex]
\Bi(x) \sim & \; B(x) - \ii \, A(x) \,, 
\nonumber\\
& \; |x| \to \infty \,, \quad 
-\pi < \arg(x) \leq -\pi + \epsilon \,.
\end{align}
\end{subequations}
The asymptotics for the 
derivative of the Airy $\Bi$ function are 
obtained by making the following replacements
in Eq.~\eqref{BiAsymp}: $\Bi \to \Bi'$,
$A \to C$, and $B \to D$.

The ``principal of asymptotic overlap'' 
which has been employed for a number of 
special functions in Ref.~\cite{CaEtAl2007}
would now call for an application of the 
power series~\eqref{AiSmall}, \eqref{ApSmall}, 
\eqref{BiSmall} and~\eqref{BpSmall}
for small $|x| < R_0$, and for the use of the 
asymptotic expansions~\eqref{AiAsymp},~\eqref{ApAsymp}
and~\eqref{BiAsymp} for $|x| > R_1$,
with the region $R_0 < |x| < R_1$ being 
bridged by the application of a nonlinear sequence
transformation to the 
asymptotic expansions~\eqref{AiAsymp},~\eqref{ApAsymp},
and~\eqref{BiAsymp}.
However, the situation is more complicated 
in reality, and it is impossible to choose
the parameters $R_0$ and $R_1$ uniformly in the 
complex plane, independent of the complex phase 
of the argument $x$.

A possible scheme which is sufficient 
to obtain at least 20 digits of accuracy for 
all $x$ is outlined in Fig.~\ref{fig3}.
The designation ``POWER'' indicates the 
use of the power series~\eqref{AiSmall}, 
\eqref{ApSmall}, \eqref{BiSmall} and~\eqref{BpSmall},
as appropriate, the designation ``ASYMP'' 
indicates the application of the 
asymptotic expansions~\eqref{AiAsymp},~\eqref{ApAsymp},
and~\eqref{BiAsymp}, where the series $A(x)$ and $B(x)$
are given in Eqs.~\eqref{defAx} and~\eqref{defBx}.
Finally, the designation ``TRAFO'' needs to be 
explained: it denotes the application of a
nonlinear sequence transformation to the 
asymptotic series, with the aim of extending their
validity to lower modulus of the argument 
than what they would otherwise be valid for.

A rather thorough discussion of
an application of sequence transformation to a 
nontrivial problem has been given in Section~2.4.2 of Ref.~\cite{CaEtAl2007},
in the context of the relativistic Green 
function for the hydrogen atom. 
Essentially, sequence transformations are generalizations
of Pad\'{e} approximations
that fulfill accuracy-through-order relations, i.e.,
they constitute rational functions that
reproduce the first few terms of a given 
input series when expanded back in powers of the 
argument. By reformulating the sequence transformation 
in terms of optimized remainder estimates~\cite{We1989},
one can enhance the rate of convergence as compared to Pad\'{e}
approximations. A more thorough discussion of sequence
transformation would go beyond the 
scope of the current article.
The sequence transformation used in the current 
article is defined as follows.
Let $s_n$ denote the $n$th partial sum of a series
\begin{equation}
s_n = \sum_{k=0}^n a_k \,,
\qquad 
m= 1,2,\dots,
\end{equation}
where the $a_n$ are the terms in the series to be summed,
and in our case, the terms in the asymptotic 
expansions~\eqref{defABCD}. Let the difference operator $\Delta$ be 
explained as
\begin{equation}
\Delta s_n = s_{n+1} - s_n = a_{n+1} \,.
\end{equation}
As described in Refs.~\cite{We1989,CaEtAl2007,We2010}
(see also \S 3.9 of Ref.~\cite{OlLoBoCl2010}),
in many cases the following sequence transformation
(Weniger transformation),
\begin{equation}
{\delta}_k^{(n)} (\beta, s_n) = 
\frac
{\displaystyle
\sum_{j=0}^{k} ( - 1)^{j} 
\left( \begin{array}{c} k \\ j \end{array} \right) 
\frac {(\beta + n +j )_{k-1}} {(\beta + n + k )_{k-1}} 
\frac {s_{n+j}} {\Delta s_{n+j}} }
{\displaystyle
\sum_{j=0}^{k} ( - 1)^{j} 
\left( \begin{array}{c} k \\ j \end{array} \right) 
\frac {(\beta + n +j )_{k-1}} {(\beta + n + k )_{k-1}} 
\frac {1} {\Delta s_{n+j}} } 
\label{dWenTr}
\end{equation}
leads to an accelerated convergence
(or summation in the case of divergence) 
of the input sequence $\{ s_n \}_{n=0}^\infty$
of partial sums of the series of the $a_k$.
Here, $\beta$ is a shift parameter
which we choose as $\beta = 1$
(as given in Refs.~\cite{We1989,CaEtAl2007}).
The starting order for the transformation 
in Eq.~\eqref{dWenTr} is $n$ which we choose to be $n=0$.
Then, ${\delta}_k^{(0)} (\beta, s_0)$
denotes the $k$th order Weniger transform
of the input sequence $\{ s_n \}_{n=0}^\infty$.
The central idea in the construction of the 
sequence transformation~\eqref{dWenTr} is the 
expansion of the remainder term in an inverse factorial
series, as explained in Ref.~\cite{We1989}.
The use of recursion relations described 
in Section~8.3 of Ref.~\cite{We1989} is imperative in order
to ensure the computational efficiency and numerical
stability in the computations of the Weniger 
transforms. In Refs.~\cite{We1989,We1996c,CaEtAl2007},
it has been established that the Weniger transformation 
is very powerful at summing the factorially
divergent series that result from the asymptotic 
expansions of special functions in terms of 
hypergeometric series.

We find, by comparison to calculations with 
extended arithmetic using computer algebra
systems~\cite{Wo1988} and by comparison of different
methods along the separating domains indicated
in Fig.~\ref{fig3}, that 
naive termination criteria are sufficient in order
to reach a prescribed accuracy of 20 decimals 
in the entire complex plane, for the Airy functions.
By a naive termination criterion, we mean that the 
summation of terms in the power series, 
or the summation of terms in the asymptotic 
series, is terminated when the next higher-order term
divided by the most recently calculated partial sum
is a factor 100 less than the prescribed 
accuracy for the evaluation of the Airy function.
Likewise, the calculation of subsequent higher-order
transforms ${\delta}_k^{(0)}$ of the 
Weniger transform~\eqref{dWenTr} is terminated
when the apparent convergence of three consecutive transform
(maximum absolute value of the 
relative difference of transforms $k$, $k+1$ and $k+2$)
is smaller than the prescribed accuracy by a factor of 100.

The appropriate evaluation method for given 
complex argument $x$ is indicated in 
Fig.~\ref{fig3}, 
as a function of $\Re(x)$ and $\Im(x)$.
The derivatives $\Ai'(x)$ and $\Bi'(x)$ are
evaluated using the same schemes as 
$\Ai(x)$ and $\Bi(x)$, respectively,
but (for the asymptotic regime of large argument)
with the replacements $A\to C$, and $B \to D$
[see Eq.~\eqref{defABCD}].
We should also clarify the exact formulas for 
some of the separating curves in Figs.~\ref{fig3}(a)
and~\ref{fig3}(b). 
The curved separating line between ``TRAFO'' 
and ``POWER'' in Fig.~\ref{fig3}(a)
follows the formula
\begin{equation}
|x| < 5 + \frac{15}{\pi} \, |\arg(x)| \,,
\quad
|\arg(x)| < \frac{2\pi}{3}  \,,
\end{equation}
and it meets the outer asymptotic
region ``ASYMP'' at $|\arg(x)| = \frac{2\pi}{3}$
and $|x| = 15$. In Fig.~\ref{fig3}(a), the transition from the 
asymptotic series $A$ to the asymptotic
series $A + \ii \, B$ takes place at 
$\arg(x) = \pm \frac{5\pi}{6}$.
In Fig.~\ref{fig3}(b), for the $\Bi$ integral,
the power series is used for $|x| < 5$ uniformly for 
any complex phase of $x$. For $|x| < 15$, the 
power series is also used, but only for 
$|\arg(x)| < \pi/3$.
The transition from the asymptotic series 
$2\, B$ to the asymptotic
series $2\, B + \ii A$ takes place for 
$|x| \geq 15$, and $\arg(x) = \pm \pi/6$.

%
%
\begin{figure}[t!]
\begin{center}
\includegraphics[width=0.55\linewidth]{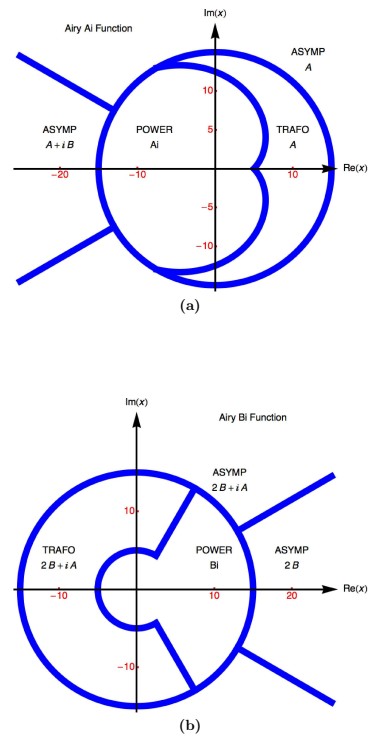} \\
\end{center}
\caption{\label{fig3}
(Color online.)
Algorithms used for the 
Airy $\Ai$ function [panel (a)] and 
for the Airy $\Bi$ function [panel(b)].
The explanation is in the text.}
\end{figure}

%
%
\subsection{Summation of uniform asymptotics}
\label{special}

The uniform asymptotic 
formulas~\eqref{Juni},~\eqref{Yuni},
\eqref{H1uni} and~\eqref{H2uni} 
are expansions for large $\nu$ for 
argument $z = \nu \, y$ of the Bessel functions.
They can be written in a form that is amenable to 
the application of the convergence acceleration
transformation~\eqref{dWenTr}, as follows.
We consider as an example Eq.~\eqref{Juni}. 
Then,
\begin{subequations}
\label{sum_form}
\begin{align}
J_\nu(\nu \, y) \sim & \; \sum_{k=0}^\infty \tilde a_k \,,
\qquad
\tilde s_n = \sum_{k=0}^n \tilde a_k \,,
\\[2ex]
\tilde a_k =& \;
\left( \frac{4 \zeta}{1 - y^2} \right)^{1/4}  \, \left\{
\frac{ \Ai(\nu^{2/3} \, \zeta) }{\nu^{1/3}} \;
\frac{ a_k(\zeta) }{ \nu^{2 k} }  \right.
\nonumber\\[1ex]
& \left. + \frac{ \Ai'(\nu^{2/3} \, \zeta) }{\nu^{5/3}} \;
\frac{ b_k(\zeta) }{ \nu^{2 k} } \right\} \,.
\end{align}
\end{subequations}
Here, the Airy functions need to be evaluated only 
once as they do not depend on $k$. Treating the Airy 
function term and the derivative term together
(i.e., as a single term $\tilde a_k$),
we empirically observe better convergence in many cases.
Written in the form~\eqref{sum_form}, the convergence of the 
uniform asymptotic for the Bessel $J$ function can be 
accelerated by calculating the 
transforms
${\delta}_k^{(n)} (\beta, \tilde s_n)$
where again, we choose $n=0$ in the current 
investigation. We empirically observe that
the use of the convergence
accelerator does not induce any stability 
issue even in cases ($\nu \neq z$) where the 
apparent convergence of the asymptotic series
would be sufficient to calculate the 
Bessel $J$ function to the required accuracy.

A special case still necessitates a modification
of the algorithm. Namely, for large and almost equal $\nu$ and $z$,
the confluence of the two saddle points defining
the Bessel function~(see Appendix~\ref{appa} below)
leads to numerically prohibitive cancellations 
in evaluating the series of the 
$a_k(\zeta)$ and $b_k(\zeta)$ coefficients.
It is impossible to overcome these difficulties
with the necessarily finite precision of computer
systems in the limit $z \to \nu$;
indeed, as described below in Section~\ref{turning},
the asymptotic expansions for large $\nu$ at
exact equality $\nu = z$ are obtained after the cancellation
of a number of divergent terms in $\epsilon$ after setting
$\nu = z + \epsilon$.

Eventually, we find that the region
in very close vicinity of the turning point $\nu = z$
can only be overcome by an explicit use of the 
recurrence relation
\begin{equation}
\label{recgen}
{\cal J}_{\nu - 1}(x) + {\cal J}_{\nu + 1}(x) = 
\frac{2 \nu}{x} \, {\cal J}_\nu(x) \,,
\end{equation}
with the notion of displacing $\nu$ from $z$,
where again ${\cal J}$ stands for $J$, $Y$, $H^{(1)}$ or $H^{(2)}$.
We find that numerical cancellations are most severe 
when $\nu$ is slightly lower than $\Re(z)$.
So, for $\Re(z) - 21 < \nu < \Re(z) + 2$, we therefore express 
${\cal J}_\nu(z)$ as a function of ${\cal J}_{\nu+24}(z)$ and
${\cal J}_{\nu+25}(z)$ by repeated use of the 
recurrence relation~\eqref{recgen}.
Indeed, by repeated application of 
Eq.~\eqref{recgen} one can express
${\cal J}_\nu(z)$ as
\begin{equation}
\label{by24}
{\cal J}_\nu(z) = 
f_{24}(\nu, z) \; {\cal J}_{\nu+24}(z) +
f_{25}(\nu, z) \; {\cal J}_{\nu+25}(z) \,,
\end{equation}
where $f_{24}$ and $f_{25}$ are somewhat lengthy 
rational functions of their arguments.
Their explicit form can easily be obtained using
computer algebra~\cite{Wo1988}.

The recurrence relation~\eqref{by24} is applied in the direction of 
increasing order $\nu$ of the Bessel function, where it
may be be unstable [indeed, 
$J_\nu(z) \to 0$ and $Y_\nu(z) \to \infty$
for $\nu \to \infty$ at constant $z$].
One might thus expect that the shift~\eqref{by24}
could induce numerical instability.
However, because the shift is applied in the region 
$\nu \approx z$, where the Bessel function is not yet
in its asymptotic regime for 
large $\nu \gg z$, the numerical cancellations are
only minor and do not constitute a matter of concern.

%
%
\subsection{Numerical examples}
\label{numerical}

\begin{figure}[t!]
\includegraphics[width=0.9\linewidth]{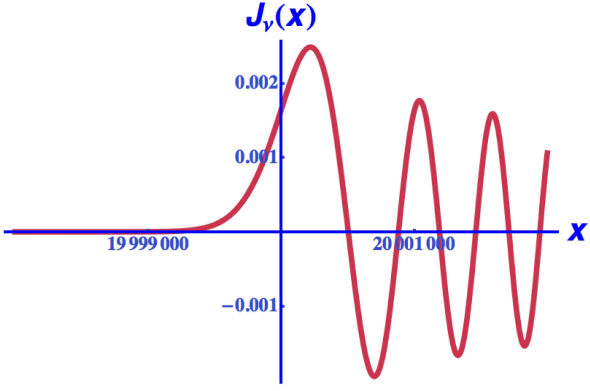}
\caption{\label{fig4} Plot of the Bessel function 
$J_\nu(x)$ with $\nu = 20,000,000.2$ in the 
turning point region.}
\end{figure}

\begin{figure}[t!]
\includegraphics[width=0.9\linewidth]{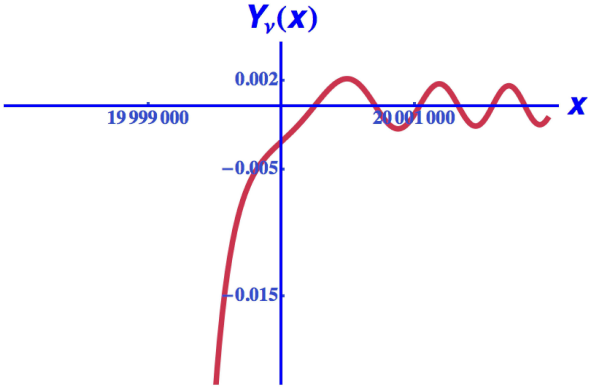}
\caption{\label{fig5} Plot of the Bessel function
$Y_\nu(x)$ with $\nu = 20,000,000.2$ in the
turning point region.}
\end{figure}

In order to demonstrate the power of 
the algorithms proposed in the current article,
we now turn our attention to a number of numerical example cases.
We first present two evaluations for non-integer order
and argument of Bessel functions, in the range $\nu \approx z$,
which result in
\begin{equation}
J_{5\,000\,000.2}(5\,000\,000.1) = 
2.614\,463\,954\,691\,926 \times 10^{-3} \,,
\end{equation}
and in
\begin{equation}
Y_{5\,000\,000.2}(5\,000\,000.1) = 
-4.533\,251\,771\,400\,041 \times 10^{-3} \,.
\end{equation}
For an example with $\arg(z) = \pi/3$, we obtain
\begin{align}
\label{eq:example_Ynu}
& H^{(1)}_{5\,000\,000.2}(5\,000\,000.1 \, \exp(\ii \pi/3) )
\nonumber\\[1ex]
& \qquad = 
-6.120\,398\,939\,598\,734 \times 10^{-954990}  
\nonumber\\[1ex]
& \qquad \qquad - \ii \, 
 1.992\,559\,471\,616\,042 \times 10^{-954989}  \,.
\end{align}
While the result is nearly zero, numerical
evaluations of this kind are needed because 
other terms in angular momentum expansions in quantum
electrodynamics lead to extremely slowly convergent 
series~\cite{JeMoSo1999,JeMoSo2001pra} due to 
an interplay of increasing and decreasing terms
(as the angular momentum quantum number is increased).
Another example of explicit evaluation of the 
Hankel function for extremely large 
order and argument in the 
turning point regime $\nu \approx z$ reads,
\begin{align}
& H^{(1)}_{6\,000\,000.2}(6\,000\,000.7) 
\\[1ex]
& \qquad = 
2.467\,848\,322\,382\,092 \times 10^{-3}
\nonumber\\[1ex]
& \qquad\qquad
- \ii \, 4.252\,887\,224\,934\,845 \times 10^{-3}\,,
\nonumber
\end{align}
and
\begin{align}
& H^{(2)}_{6\,000\,000.2}(6\,000\,000.7) = 
\\[1ex]
& \qquad 
2.467\,848\,322\,382\,092 \times 10^{-3}
\nonumber\\[1ex]
& \qquad\qquad
+ \ii \, 
4.252\,887\,224\,934\,845\times 10^{-3} \,.
\nonumber
\end{align}
For $r \approx 1$, it is instructive to 
verify the Bessel functions using the following sum rule,
\begin{equation}
\label{sum_rule}
\frac{\exp\bigl( -y[1-r] \bigr)}{y[1-r]} \; = \; - \,
\sum_{\ell=0}^{\infty} (2l+1) \, j_\ell (\ii \, r \, y) \, 
h^{(1)}_\ell (\ii \, y) \,,
\end{equation}
for $r \in (0,1)$ and $y > 0$.
The sum over $\ell$ constitutes a slowly convergent series whose 
can otherwise by accelerated using the 
combined nonlinear-condensation transformation~\cite{JeMoSoWe1999}.
The sum rule directly follows from the angular
momentum expansion of the Green function of the Helmholtz
equation for the case of collinear arguments,
as given in Chap.~9 of Ref.~\cite{Ja1998}.
It can also be derived as a reformulation of 
Eq.~(10.60.3) of Ref.~\cite{OlLoBoCl2010}.
We have checked our values of the Bessel functions on the basis 
of this sum rule and plotted Bessel functions of high order
in the turning point region (see Figs.~\ref{fig4} and~\ref{fig5}).

%
%
\section{ASYMPTOTICS AT THE TURNING POINT}
\label{turning}

The direct application of the uniform asymptotics
becomes problematic when the argument and the order of a
Bessel function are almost equal,
because of numerical cancellations involved in 
evaluating the individual coefficients in this 
case. We remember that for the case $\nu \approx z$, 
special methods have to be employed for the 
summation of the uniform asymptotic 
expansion~(see Section~\ref{special}).
This is because of highly significant numerical
cancellations in the evaluation of the 
sums defining $a_k$ and $b_k$ given in Eqs.~\eqref{ak}
and~\eqref{bk}. These numerical cancellations reflect
the confluence of the two saddle discussed below
in Appendix~\ref{appa} and are present even if the 
sums in Eqs.~\eqref{ak} and~\eqref{bk} are finite.

Nevertheless it is possible to investigate the limit $\nu \to z$ by an analytic
expansion setting $\nu = z + \epsilon$. In the calculation, the result is
obtained only after the cancellation of a few divergent terms which are
inverse powers of $\epsilon$.  For exact equality $\nu = z$, it is possible to
derive analytic approximations. These serve as an important test device for
the numerical algorithms.  The results given in formulas (9.3.31)--(9.3.34) of
Ref.~\cite{AbSt1972} for lower order terms in the asymptotic expansions do not
reach a sufficiently high order in order to be useful for the current
investigation, and we have thus calculated a few more terms in the asymptotic
expansions of $J_\nu(\nu)$ and $Y_\nu(\nu)$ for $\nu \to \infty$. 

The general structure of the asymptotic expansions is as follows,
\begin{subequations}
\begin{align}
J_\nu(\nu) =& \; 
\frac{a}{\nu^{1/3}} \, 
\sum_{k=0}^{\infty} \frac{\alpha_k}{\nu^{2 \, k}} 
- \frac{b}{\nu^{5/3}} \, 
\sum_{k=0}^{\infty} \frac{\beta_k}{\nu^{2 \, k}} \,,
\\[1ex]
Y_\nu(\nu) =& \;
-\frac{3^{1/2} \, a}{\nu^{1/3}} \,
\sum_{k=0}^{\infty} \frac{\alpha_k}{\nu^{2 \, k}}
- \frac{3^{1/2} \, b}{\nu^{5/3}} \,
\sum_{k=0}^{\infty} \frac{\beta_k}{\nu^{2 \, k}} \,,
\\[1ex]
J'_\nu(\nu) =& \;
\frac{b}{\nu^{2/3}} \,
\sum_{k=0}^{\infty} \frac{\gamma_k}{\nu^{2 \, k}}
- \frac{a}{\nu^{4/3}} \,
\sum_{k=0}^{\infty} \frac{\delta_k}{\nu^{2 \, k}} \,,
\\[1ex]
Y'_\nu(\nu) =& \;
= \frac{3^{1/2} \, b}{\nu^{2/3}} \,
\sum_{k=0}^{\infty} \frac{\gamma_k}{\nu^{2 \, k}}
+ \frac{3^{1/2} \, a}{\nu^{4/3}} \,
\sum_{k=0}^{\infty} \frac{\delta_k}{\nu^{2 \, k}} \,.
\end{align}
\end{subequations}
The constants $a$ and $b$ are given by~\cite{AbSt1972}
\begin{equation}
a = \frac{2^{1/3}}{3^{2/3} \, \Gamma(2/3)} \,,
\end{equation}
and
\begin{equation}
b = -\frac{2^{2/3}}{3^{1/3} \, \Gamma(1/3)} \,.
\end{equation}
Numerical values of the coefficients 
$\alpha_k$, $\beta_k$, $\gamma_k$, $\delta_k$
were given in (unnumbered) equations following 
Eq.~(9.3.34) of Ref.~\cite{AbSt1972}, but  only 
in numerical form. 
Using the formalism outlined in Refs.~\cite{AbSt1972,Wa1922,Ma1993bessel},
and computer algebra~\cite{Wo1988}, it is relatively 
easy to evaluate these coefficients analytically.
The first terms read
\begin{align}
\alpha_0 =& \; 1 \,, \qquad
\alpha_1 = -\frac{1}{225} \,, \qquad
\alpha_2 = \frac{151439}{218295000} \,,
\nonumber\\[1ex]
\alpha_3 =& \; -\frac{887278009}{2504935125000} \,, 
\nonumber\\[1ex]
\alpha_4 =& \; \frac{1374085664813273149}{3633280647121125000000} \,,
\nonumber\\[1ex]
\alpha_5 =& \; -\frac{1065024810026227256263721}%
{1540965154460247140625000000} \,.
\end{align}
for the $\alpha_k$ coefficients. 
For $\alpha_j$ with $j = 6,7,8$, the integers in the 
numerator and denominator are too long to be
displayed, practically. The results read, up to 
30~decimals,
\begin{align}
\alpha_6 =& \; 
 0.00192\,82196\,42487\,75701\,38042\,30112 \,,
\nonumber\\[1ex]
\alpha_7 =& \; 
-0.00762\,60912\,66562\,73551\,11507\,07644 \,,
\nonumber\\[1ex]
\alpha_8 =& \; 
 0.04059\,16252\,02439\,02610\,46110\,96353 \,.
\end{align}
We obtain for the $\beta_k$,
\begin{align}
\beta_0 =& \; \frac{1}{70} \,, \qquad
\beta_1 = -\frac{1213}{1023750} \,, 
\nonumber\\[1ex]
\beta_2 =& \; \frac{16542537833}{37743205500000} \,,
\beta_3 = -\frac{9597171184603}{25476663712500000} \,, 
\nonumber\\[1ex]
\beta_4 =& \; \frac{53299328587804322691259}{91182706744837207500000000 } \,,
\nonumber\\[1ex]
\beta_5 =& \; -\frac{70563256104582737796094772987}%
{49341242187300033908437500000000} \,.
\end{align}
For $\beta_j$ with $j = 6,7,8$, the results read, up to
30~decimals,
\begin{align}
\beta_6 =& \; 
 0.00506\,84595\,77410\,25775\,49192\,90289 \,,
\nonumber\\[1ex]
\beta_7 =& \; 
-0.02455\,31387\,44039\,61140\,77960\,53629 \,,
\nonumber\\[1ex]
\beta_8 =& \; 
 0.15584\,22313\,83604\,27406\,01873\,20654\,.
\end{align}
The results for the $\gamma_k$ read,
\begin{align}
\gamma_0 =& \; 1 \,, \qquad
\gamma_1 = \frac{23}{3150} \,, \qquad
\gamma_2 = - \frac{604523}{644962500} \,,
\nonumber\\[1ex]
\gamma_3 =& \; \frac{2264850139339}{5095332742500000} \,, 
\nonumber\\[1ex]
\gamma_4 =& \; -\frac{160913976870912403}{353106559055250000000} \,,
\nonumber\\[1ex]
\gamma_5 =& \; \frac{216363828773939104866579281}%
{266709417228648831937500000000} \,.
\end{align}
For $\gamma_6$, $\gamma_7$, and $\gamma_8$, we again
give numerical results, up to 30~decimals,
\begin{align}
\gamma_6 =& \; -0.00222\,15829\,52781\,90513\,99808\,25549\,,
\nonumber\\[1ex]
\gamma_7 =& \;  0.00866\,43390\,74500\,57903\,06978\,05128\,,
\\[1ex]
\gamma_8 =& \; -0.04561\,48821\,43058\,95740\,37593\,20268\,.
\nonumber
\end{align}
Finally, we obtain for the $\delta_k$,
\begin{align}
\delta_0 =& \; \frac{1}{5} \,, \qquad
\delta_1 = -\frac{947}{346500} \,, 
\\[1ex]
\delta_2 =& \; \frac{11192989}{18555075000} \,,
\delta_3 = -\frac{100443412440047}{262141460831250000} \,, 
\nonumber\\[1ex]
\delta_4 =& \; \frac{11007971229145235539}{22905464949241875000000} \,,
\nonumber\\[1ex]
\delta_5 =& \; -\frac{180026595127347603856424798473}%
{180354561678027325338750000000000} \,.
\nonumber
\end{align}
The first 30 decimals of the results for $\delta_{6,7,8}$ read
\begin{align}
\delta_6 =& \;  0.00309\,74954\,54537\,99792\,65245\,46715\,,
\nonumber\\[1ex]
\delta_7 =& \; -0.01341\,95416\,64085\,43255\,61825\,95640\,,
\nonumber\\[1ex]
\delta_8 =& \;  0.07735\,88016\,72766\,55151\,01508\,73624 \,.
\end{align}
We have used these analytic formulas in our 
verification of numerical results obtained using the 
method described in Section~\ref{summation}.

%
%
\section{CONCLUSIONS}
\label{conclu}

We have described an algorithm for the evaluation of an individual Bessel or
Hankel function for large order $\nu$, and arbitrary complex argument $z$.  The
method relies on the use of the uniform asymptotic expansions given in
Eqs.~\eqref{Juni},~\eqref{Yuni}, \eqref{H1uni}~and~\eqref{H2uni}.  These
asymptotic expansions involve $\Ai$ and $\Bi$ functions, and their
derivatives.  Consequently, in passing and also as a prerequisite for our
algorithm, we need to develop numerical code for the Airy functions valid in
the entire complex plane.  This is described in Section~\ref{airy}, and the
appropriate algorithms for the different regions in the complex plane are given
in Fig.~\ref{fig3}.  Essentially, these represent an adapted
implementation of the ``principle of asymptotic overlap'' where a power series
is used for small argument, an asymptotic expansion is used for large argument,
and the region in between is bridged by a nonlinear sequence transformation
[Weniger transformation, see Eq.~\eqref{dWenTr}]. 

As evident from
Fig.~\ref{fig3}, this principle needs to be adapted for the
Airy functions, and the regions separating the algorithms depend on the complex
phase of the argument. Using the Airy function algorithm and the uniform
asymptotic expansion, we obtain expressions for the Bessel and Hankel functions
which can readily be evaluated for almost arbitrary values of $\nu$ and $z$.
However, formidable numerical cancellations still prevent us from using the
uniform asymptotic expansions directly when $\nu \approx z$. In this case, the
region near $\nu = z$ is bridged by the application of a convergence
accelerator (Weniger transformation), applied in this case to the infinite
asymptotic series defining the expansion for large $\nu$ in the uniform
asymptotic formulas~\eqref{Juni},~\eqref{Yuni},
\eqref{H1uni}~and~\eqref{H2uni}. The transformation is applied 
after the order $\nu$ is shifted away [see Eq.~\eqref{by24}] 
by a finite amount from the 
argument $z$ using the recurrence formula satisfied by the Bessel
function~\eqref{recur}.
Numerical examples are given in Section~\ref{numerical}.

In a typical application within physics, one needs Bessel functions for both
real as well as complex argument.  For example, in bound-state quantum
electrodynamics (QED), one needs to describe complex photon energies.
Therefore, it is highly desirable to have an algorithm that is applicable in
the entire complex plane, for a Bessel function of large order. Here, we use a
combination of ideas to overcome the severe difficulties faced by this
endeavor.  We use a nonlinear sequence transformation for the summation of
asymptotic expansions for the Airy functions in order be able to use this
expansion for small and moderate argument where the usual paradigm of
truncating the asymptotics at the smallest term would otherwise yield
dissatisfactory accuracy.  The Stokes phenomenon must still be taken into
account in terms of the complex phase of the argument of the Airy function.
Finally, we overcome the remaining numerical difficulties in the summation of
the uniform asymptotic expansion, associated with the turning point $\nu = z$,
by a simple shift of the order versus the argument, and by the further use of
the Weniger transformation for the summation of the uniform asymptotic
expansion in the transition region near $\nu = z$.

Examples where Bessel functions of large argument and order are needed, include
scattering problems and calculations with photon propagators and fermion
propagators in atomic physics and field theory, and diverse technical
application areas.  Although it is possible to recursively evaluate arrays of
Bessel functions of the same argument with varying order, the individual
evaluation of Bessel functions in ``extreme'' argument ranges remains an
elusive problem.  Here, we attempt to address this problem via a combination of
summation algorithms, convergence accelerators, recurrence relations and the
``principle of asymptotic overlap,'' adapted to the problem at hand.  In
passing, we also address the numerically accurate calculation of Airy functions
in the entire complex plane.
We leave it to the interested reader to develop their own implementation
of the methods described here, adapted to the arithmetic accuracy
requirements for the particular application in question.

Finally, let us indicate a few open problems in the area.
We have outlined the general algorithm used in our 
evaluation of Bessel functions. The description of an implementation,
including an example code, possibly with more refined 
termination criteria than those indicated in 
Section~\ref{airy} for the Airy functions, would certainly be
of value to the scientific community. 
Secondly, in view of the considerations outlined below in 
Appendix~\ref{appa}, it seems feasible to 
construct an equally powerful algorithm based on an 
adaptation of the saddle point integration.
In this case, considerable care needs to be vested into the 
turning point case $\nu \approx z$ as well.
Finally, preliminary investigations (not described
here in any further detail) indicate that
the Debye expansion, which is different from the 
uniform asymptotic expansion, may be used to good effect
for the case of real $\nu$, and real argument $z$,
of the Bessel function. The Debye expansion is given in
Eqs.~(9.3.7),.~(9.3.8),~(9.3.11) and~(9.3.12)
of Ref.~\cite{AbSt1972}. The Debye expansion may be 
combined with convergence accelerators.
It would be instructive and worthwhile to construct
a complementary algorithm for Bessel and Hankel
functions valid only for real argument.

%
%
\section*{Acknowledgments}

U.D.J.~acknowledges support by the National Science Foundation
and support by a precision measurement grant
from the National Institute of Standards and Technology.
The work of E.~L.~was supported by  the Japan Society for the Promotion of
Science (Grant-in-Aid for Scientific Research No. 21-09238).
U.D.J.~acknowledges the kind hospitality of the 
National Institute of Standards and Technology during 
the month of August 2010, where part of this work was completed.

\appendix

%
%
\section{SADDLE POINT INTEGRATIONS}
\label{appa}

%
%
\subsection{Orientation}

In this Appendix, we explore alternative numerical procedures, with a special
emphasis on numerical integration around saddle points.  Before we come to a
discussion of the procedure involved, let us briefly review other, alternative
approaches to the calculation of Bessel and Hankel functions which have been
discussed in the literature.  In Ref.~\cite{Te1994}, it has been advocated to
directly integrate the defining differential equation in suitable directions in
the complex plane.  In Ref.~\cite{ScAnGo1979}, an analogous approach has been
advocated for the calculation of Airy functions.  There is a further obvious
problem at the turning point, where the argument $z$ and order $\nu$ of the
Bessel function are nearly equal. In this case, there are huge numerical
cancellations in the calculation of the expansion coefficients that multiply
the Airy functions which are required in order to evaluate the uniform
asymptotic expansions.  In Ref.~\cite{Te1997uniform}, it has been suggested to avoid
using the expansion in inverse powers of the order of the Bessel function in
this case, and to use a convergent expansions in terms of a scaled complex
argument $\zeta$ of the special function.  The procedure advocated in
Ref.~\cite{Te1997uniform} is applicable to the transition region $\nu \approx z$ where
$\zeta \approx 0$, but cannot be universally applied when $\zeta$ is manifestly
different from zero.  Another method~\cite{HoOD1999} is to expand the
coefficients of the asymptotic expansion into hyperasymptotics. Yet, the
asymptotic expansions of the coefficients of the asymptotic expansions (sic!)
are not applicable to all cases of interest and introduce another level of
complexity into the problem.  Hadamard series expansions as an alternative to
the uniform asymptotics (still equal in limiting cases) have been investigated
in Refs.~\cite{Pa2001hadamard1,Pa2001hadamard2,Pa2004hadamard3}.
Finally, let us also mention that some special parameter 
cases of interest
have received attention in the literature. For example,
asymptotic expansions for Bessel functions of the third kind
of imaginary order have been discussed in
Refs.~\cite{Ba1966,ShWo2009}.

However, the preeminent issue in the formulation of 
possible alternative evaluation methods for 
Bessel and Hankel functions seems to be connected
with a numerical integration about saddle points in the 
complex plane. Expansions about the saddle points 
give rise to the above mentioned 
uniform asymptotic expansions~\eqref{Juni},~\eqref{Yuni},
\eqref{H1uni} and~\eqref{H2uni}. In general, the
evaluation of Hankel, and Bessel functions can be written in terms of two
saddle points which dominate the paths of integration.  Expanding to second
order about the saddle point, we obtain linear contours that approximate the
paths of steepest descent and may be used in order to evaluate the functions,
approximately.  This approach, which has been advocated previously in
\cite{Ma1993bessel} and recently in \cite{SmdHTi2009}, is appealing because of its
relative simplicity.
Since the integrand of the integral defining $J_\nu (z)$ is
in general highly oscillating, the integration paths
must be chosen with care. Here, we describe a somewhat simplified
approach, in which one approximates the ideal path by a short linear segment
crossing the saddle point. 
Our treatment is similar to that
in \cite{SmdHTi2009}, however, we pay special attention to the case when
argument and order of the Bessel function are almost equal, which was not
considered in \cite{SmdHTi2009}.

In the case of equal argument and order, $\nu = z$, the phenomenon of
confluence of the two saddle points needs to be analyzed.  For $\nu \approx z$,
the two saddle points approach each other.  For precise equality $\nu = z$, the
two saddle points coalesce, and there are three paths of steepest descent,
complemented by three paths of steepest ascent, about the saddle point. This
leads to a more complicated situation with kinks in the contours of steepest
descent, which must be used in order to evaluate the special functions.  
The investigations reported below are therefore of more general 
interest with respect to the numerical problems related to the 
confluence of the saddle points.

%
%
\subsection{Complex integral representations}

The defining integrals for the Hankel functions $H^{(1)}_\nu (z)$ and
$H^{(2)}_\nu (z)$ read  ($\real z>0$, and arbitrary complex $\nu$)
\begin{subequations}
\begin{align}
\label{eq:Hankeldef}
H^{(1)}_\nu (z) =& \; 
\frac{1}{\pi\ii} \, 
\int_{-\infty-\ii\pi}^{\infty} \exp\left( g(t) \right) \dd t,
\\
\label{eq:Hankeldef2}
H^{(2)}_\nu (z)=& \;
\frac{1}{\pi\ii}\int^{-\infty+\ii\pi}_{\infty}
\exp\left( g(t) \right) \dd t,
\end{align}
\end{subequations}
with
\begin{equation}
\label{eq:gdef}
g(t) = -z\sinh(t)+\nu t = -\nu \, (y\sinh(t) - t) \, .
\end{equation}
Here we have introduced the variable $y=z/\nu$, 
which will prove convenient in
the ensuing discussion.  Combining the two Hankel functions to the Bessel
function $J_\nu (z)=[H^{(1)}_\nu (z)+H^{(2)}_\nu (z)]/2$, we obtain from
Eq.~\eqref{eq:Hankeldef} the integral representation for $J_\nu (z)$:
\begin{equation}
\label{eq:Jnudef}
J_\nu (z) = \frac{1}{2\pi\ii}
\int^{-\infty + \ii \pi}_{-\infty - \ii \pi}
\exp\left( g(t) \right) \, \dd t.
\end{equation}
The contour of integration \eqref{eq:Jnudef} is arbitrary.  It can be deformed,
nevertheless, to pass through the saddle point.  Finally, the Bessel function
of the second kind is obtained as 
$Y_\nu (z)=[H^{(1)}_\nu (z)-H^{(2)}_\nu (z)]/(2\ii) =
-\ii \, [H^{(1)}_\nu (z)-J_\nu (z)]$.

The definition \eqref{eq:Hankeldef} is valid provided $\real z>0$, to ensure
convergence of the integrals. 
The identity $\sinh(t - \ii \pi) = \sinh(t + \ii \pi) = -\sinh(t)$ 
is instrumental in deriving this condition.
For $\real(z) < 0$, we apply the conversion
formulas~\eqref{six}. The case of purely imaginary $z = x + \ii \, y$ 
requires another integral representation. We may use
\begin{align}
\label{eq:modifiedIK1}
J_\nu (\ii |y|) =& \; \exp\left( \frac{1}{2}\nu \pi \ii \right) \, I_\nu (|y|), \\
\label{eq:modifiedIK2}
J_\nu (-\ii |y|) =& \;
\exp\left(-\frac{3}{2}\nu \pi \ii\right) \, I_\nu (-|y|),\\
\label{eq:modifiedIK3}
H_\nu^{(1)}(\ii |y|) =&\frac{2}{\pi \ii} \ee^{-\frac{1}{2}\nu\pi\ii}K_\nu (|y|),\\
\label{eq:modifiedIK4}
H_\nu^{(2)}(-\ii |y|) =&-\frac{2}{\pi \ii} \ee^{\frac{1}{2}\nu\pi\ii}K_\nu (|y|),
\end{align}
which expresses the Bessel functions 
in terms of modified Bessel functions $I_\nu (x)$ and $K_\nu (x)$. 
For $\real(x) >0$, these have integral representations
\begin{subequations}
\begin{align}
I_\nu (x)=& \; \frac{1}{\pi} 
\int_0^\pi \exp\left( x\cos\theta \right) \, \cos(\nu\theta) \dd\theta \nonumber\\
& \; -\frac{\sin(\nu\pi)}{\pi} \,
\int_0^\infty \exp\left( -x\cosh t -\nu t \right) \, \dd t,\\ 
K_\nu (x)=&\int_0^\infty \exp\left( -x\cosh t \right) \, \cosh(\nu t) \dd t \,.
\end{align}
\end{subequations}
Negative arguments can be transformed to positive arguments
via the relations
\begin{subequations}
\begin{align}
I_\nu (-x) =& \; \exp\left( \nu \pi \ii \right) \, I_\nu (x),\\
K_\nu (-x) =& \; \exp\left( -\nu \pi \ii \right) \, K_\nu (x)
-\pi\ii I_\nu (x).
\end{align}
\end{subequations}

%
%
\subsection{Paths of steepest descent}

We recall the definition of the integrand in the integral
representation of the Bessel function,
\begin{equation}
g(t) = -\nu \, (y \, \sinh(t) - t) \, .
\end{equation}
Expanded to second order about the saddle point $t = t_j$, 
to be specified below, one obtains
\begin{align}
\label{third}
g(t) =& \; g(t_j) + g'(t_j) \, (t-t_j) + \frac12 g''(t_j) \, (t-t_j)^2 + \ldots 
\nonumber\\
=& \; -\nu \, \left(y\sinh(t_j) - t_j\right) + 
\nu \left(1 - y \cosh(t_j) \right) \, (t-t_j) 
\nonumber\\
& \; - \frac12 \, \nu \, y \, \sinh(t_j) \, (t-t_j)^2
\nonumber\\
& \; - \frac16 \, \nu \, y \, \cosh(t_j) \, (t-t_j)^3 + \dots \,.
\end{align}
The method of numerical evaluation of the functions $J_\nu (z)$ and $Y_\nu (z)$
amounts to finding approximations to the 
integrals \eqref{eq:Hankeldef}, \eqref{eq:Hankeldef2} by
quadrature.  The starting point is to find the  saddle points $t_\pm$  of the
integrand in \eqref{eq:Jnudef}. They satisfy $g'(t_\pm)=0=y\cosh t_\pm-1$, and
are given as
\begin{subequations}
\label{tplustminus}
\begin{align}
\label{tplus}
t_+ =& \; \ln\left(\frac{1+\sqrt{1-y^2}}{y}\right)=a_++\ii \; b_+,
\\
\label{tminus}
t_- =& \; \ln\left(\frac{1-\sqrt{1-y^2}}{y}\right)=a_-+\ii  \; b_- = 
- t_+ \,,
\end{align}
\end{subequations}
for $|\imag t_\pm| \le \pi$, and the notation $t_\pm = a_\pm + \ii \, b_\pm$
(real $a_\pm$, $b_\pm$) is introduced for later use. Multiplying the 
arguments of the logarithms in Eq.~\eqref{tplustminus}, 
we see that the saddle points fulfill the relation $t_+ + t_- = 0$.
Furthermore, we have $\real t_-\le 0$ (and $\real t_+ \geq 0$).  
For the case of equal argument and order,
the two saddle point coalesce at the origin,
\begin{equation}
\nu = z \,, \qquad
t_+ = t_- = 0 \,. \qquad
\end{equation}
The path of
steepest descent (PSD) is a path $t$ in the complex plane, passing through the
saddle point, such that the imaginary part of the exponent $g(t)$  is constant
along this path. This means that the integrand 
$\exp\left( g(t) \right)$ does not oscillate.
Of course, there are further saddle points in the complex plane, described by
the formula $t_\pm+2n\pi\ii$, $n$ integer, but these are not needed in our
investigation. 

In addition, the amplitude of the integrand $\exp\left( g(t) \right)$ decreases
monotonically as we go along the path away from the saddle point.  Together,
these properties make numerical quadrature along the PSD easy. In order find an
expression for the PSD, we write 
\begin{equation}
t =\xi+\ii \, \eta \,,
\end{equation}
and require that the imaginary part of $g(t)$ 
remains constant along the path of integration.
We write $z = \zr + \ii \, \zi$, with $\zr=\real z$ and $\zi=\imag z$,
and in complete analogy, $\nu=\nur+\ii \, \nui$. 
Then,
\begin{align}
\label{eq:PSD}
& \imag g(t_\pm) = \imag g(t) \\[1ex]
& \;\; =-\zi \, \sinh\xi\, \cos\eta -
\zr \, \cosh\xi \, \sin\eta+\nui\, \xi+\nur \, \eta.
\nonumber
\end{align}
In general, Eq.~\eqref{eq:PSD} cannot
be solved analytically. The limiting behavior as $|\xi|\to\infty$
may however be deduced. Namely, for 
$|\xi|\to\infty$, one has $|\sinh \xi| \gg |\xi|$, $|\cosh \xi| \gg |\xi|$,
and $\tanh\xi \to 1$ for $\xi \to \pm \infty$.
Therefore, the asymptotic solutions are given by
\begin{subequations}
\begin{align}
\tan\eta =& \; \frac{\zi}{\zr}=\tan\left(\arg z\right),\qquad \xi\to-\infty,\\
\tan\eta =& \; -\frac{\zi}{\zr}=-\tan\left(\arg z \right),\qquad \xi\to\infty.
\end{align}
\end{subequations}
A possible approach now is to solve Eq.~\eqref{eq:PSD} numerically, and
subsequently integrate along the numerically obtained PSD \cite{SmdHTi2009}.
A helpful discussion of numerical aspects related to 
contour integrals in the complex plane is given in Chapter~5
Ref.~\cite{GiSeTe2007}.

As we shall see, to obtain modest accuracy it is enough to integrate
along a short linear segment $\Gamma$, which approximates the true PSD close to
the saddle point. 
We have previously defined the saddle points $t_j$ ($j = \pm$)
and their real and imaginary parts $a_j$ and $b_j$.
Linear approximations to the
paths of steepest descent in the complex $t$ plane
are obtained as follows,
\begin{subequations}
\label{eq:deftauofeta}
\begin{align}
t = t_j(\eta) =& \; (\eta-b_j) \; K_j + a_j + \ii \, \eta, 
\quad j = \pm \,,
\\[1ex]
t_j \equiv t_j(b_j) =& \; a_j + \ii \, b_j \,,
\qquad
\frac{\dd t_j(\eta)}{\dd \eta} = K_j + \ii \,,
\end{align}
\end{subequations}
where $K_j$ is real.
That is, we parameterize the integral in \eqref{eq:Jnudef}
(similarly for the $Y_\nu (z)$ function) 
as a function of $\eta=\imag t$ as 
\begin{align}
\label{eq:approxJnuint}
J_\nu (z) \approx & \;
\sum_{j = \pm} \frac{1}{2\pi\ii}
\int_{\Gamma_j} \exp\left( g(t) \right) \, \dd t 
\nonumber\\
=& \;
\sum_{j=\pm} \frac{1}{2\pi\ii} \int_{\eta_{\textrm{min}}}^{\eta_{\textrm{max}}}
\exp\left( g(t_j(\eta)) \right) \;
\frac{\dd t_j(\eta)}{\dd\eta} \; \dd \eta,
\end{align}
where the $\Gamma_j$ are line segments along the paths that 
define linear approximations to the contours of steepest descent,
with $j\in\{+,-\}$. The values of the 
parameters $\eta_{\textrm{min}}$ and $\eta_{\textrm{max}}$
are chosen so that a specific final accuracy is reached for the 
approximation to the integral.

The sum in \eqref{eq:approxJnuint} may run over one, or both of the saddle
points $t_\pm$, depending on the saddle point configuration.  Our expansion to
second order about the saddle point ensures that no oscillations occur in the
integrand up to this order.  Therefore, no imaginary part will be incurred in
the integrand upon expansion about the saddle point to second order
in $t - t_j$.  However, if one goes beyond second order,
oscillations will occur upon using the linear approximation~\eqref{eq:deftauofeta}.

The coefficients $K_j$ in
\eqref{eq:deftauofeta} can be found by expanding Eq.~\eqref{eq:PSD} to second
order around the saddle point $t_j$, which yields
\begin{align}
0&=\imag\left[z(t-t_j)^2\sinh(t_j)\right]\nonumber\\
&=\imag \big\{ z\left[(\xi-a_j)^2-(\eta-b_j)^2+2\ii (\xi-a_j)(\eta-b_j)\right]
\nonumber\\
&\quad\,\times [\sinh a_j \cos b_j +\ii\cosh a_j \sin b_j]\big\}. 
\end{align}
Solving for $\xi$ as a function of $\eta$, we obtain
\begin{equation}
\label{eq:xiofeta}
\xi-a_j = (\eta-b_j) \, \left(-B_j\pm\sqrt{B_j^2+1}\right) \,,
\end{equation}
where
\begin{align}
B_j=& \; \frac{\zr\sinh a_j \cos b_j - \zi\cosh a_j \sin b_j}%
{\zr\cosh a_j \sin b_j +\zi\sinh a_j \cos b_j}
\nonumber\\
=& \; \cot \arg(z\sinh t_j).
\end{align}
It will be discussed below
how to estimate the upper and lower integration limits $\eta_{\textrm{min}}$ and
$\eta_{\textrm{max}}$. The coefficient $K_j$ is given by
\begin{equation}
K_j = -B_j\pm\sqrt{B_j^2+1} \,,
\end{equation}
and the $\pm$ sign needs to be fixed.  Indeed, there are two solutions for each
saddle point, one of which corresponds to the PSD. The other solution yields
the path of steepest ascent, along which the integrand $\exp\left( g(t)
\right)$ increases. The condition 
\begin{subequations}
\begin{align}
\Psi =& \; \arg(z\sinh t_j) \in ( - \pi, \pi) \,, 
\\[1ex]
\left( K_j + \ii \right)^2 \mathop{=}^{!} & \; 
\sqrt{ K_j^2 + 1} \,
\exp\left(-\ii \, \Psi \right) \,,
\end{align}
\end{subequations}
ensures that the term of order $(t - t_j)^2 \propto \eta^2$
in Eq.~\eqref{third} leads to an exponential decrease 
(instead of increase) of the integrand around the saddle point.
It is fulfilled if we choose
\begin{subequations}
\begin{align}
K_j =& \; -B_j - \sqrt{B_j^2+1} \\
=& \; -\sqrt{ \cot^2\left( \tfrac12 \Psi \right) + 1 } 
\exp\left( - \tfrac{\ii}{2} \Psi \right) \,, \quad \Psi > 0 \,, 
\nonumber\\[1ex]
K_j =& \; -B_j + \sqrt{B_j^2+1} \\
=& \; \sqrt{ \cot^2\left( \tfrac12 \Psi \right) + 1 } 
\exp\left( - \tfrac{\ii}{2} \Psi \right) \,, \quad \Psi < 0 \,,
\nonumber
\end{align}
\end{subequations}
depending on the sign of $\Psi$.
 
\begin{figure}[t!]
\includegraphics[width=0.77\linewidth]{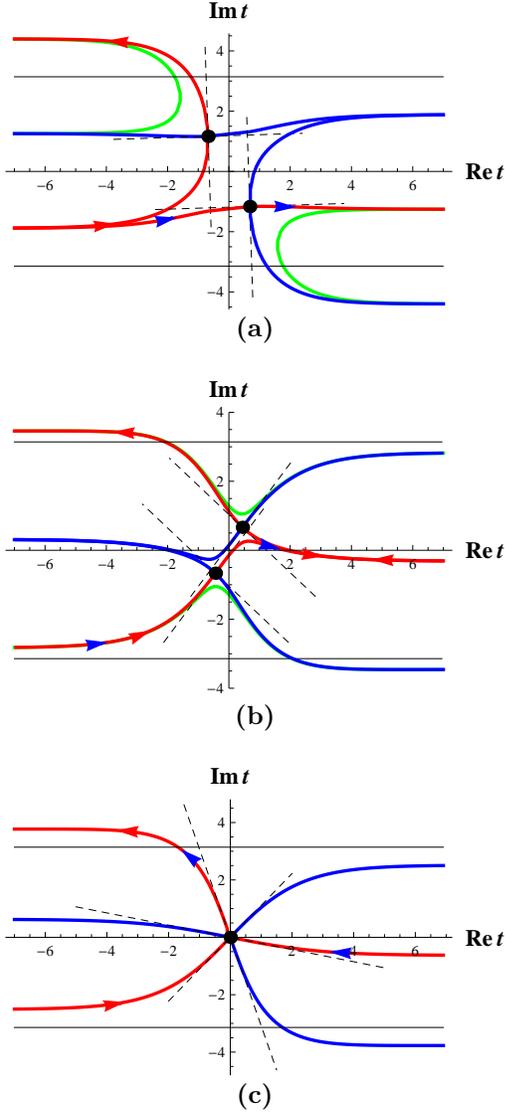}
\caption{\label{fig6} (Color online.) 
PSD for different values of the argument $z$ and
order $\nu$. In (a), $\nu=(1000+1/5) \ee^{\ii\pi/10}$, 
$z=(1200+6/25) \ee^{2\ii\pi/5}$,
in (b), $\nu=(1000+1/5) \ee^{\ii\pi/5}$, $z=(1100+11/50) \ee^{\ii\pi/10}$, and in (c),
$z=\nu=(1000+1/5) \ee^{\ii\pi/5}$. Red lines show PSDs passing through the saddle
points (black dots), blue lines paths of steepest ascent passing through the
saddle points. For orientation purposes, thin black lines are drawn at 
$\imag t=\pm\pi$. Red arrows indicate the integration path to obtain $J_\nu(z)$, and
blue arrows the paths for $H^{(1)}(z)$ [(a) and (b)] and  $H^{(2)}(z)$ [(c)].
In (a) and (b), the green curves are curves satisfying Eq.~\eqref{eq:PSD}, but
not passing through the saddle point. The slopes of the linear segments
\eqref{eq:xiofeta}, \eqref{eq:PSD_nueqz} used in the numerical evaluation are
shown with dashed lines. To obtain 10 significant figures, one needs to perform
the integration along a linear segment of length approximately equal to the
size of the black dots.}
\end{figure}

The case $\nu=z$ needs to be considered separately.
The saddle point coalesce at $t_+=t_-=0$.
The second-order terms vanishes ($\sinh t_j = 0$), and the linear approximation
to the PSD is obtained by expanding \eqref{eq:PSD} to third order.
In this case, the saddle point becomes a ``triple saddle point''
which instead of a cross formation has the shape
of a double saddle point, or of a 
``six-fold star'' with three paths of steepest descent, and 
three paths of steepest ascent in between.
The corresponding term from Eq.~\eqref{third} in this case is 
\begin{equation}
- \frac{\nu \, y}{6} \, \cosh(t_j = 0) \, (t-t_j)^3 =
- \frac{z}{6} \, t^3 \,,
\end{equation}
because,
for the confluence of the two saddle points, we have
$\nu = z$ and $t_j = 0$.  The saddle point equation thus simplifies and
reads
\begin{equation}
\imag(z \, t^3)=0.
\end{equation}
This has three solutions,
\begin{equation}
\label{eq:PSD_nueqz}
t_k=\eta(\cot\varphi_k+\ii)= \frac{\eta}{\sin\varphi_k} \, \exp(\ii\varphi_k),
\end{equation}
where $\varphi_k=(k \, \pi-\arg z)/3$, and $k\in\{0,1,2\}$, and the last
equality indicates that the confluent PSDs can be parameterized conveniently in
terms of $s = \eta/\sin\varphi_k$. The PSDs derived in Eq.~\eqref{eq:PSD_nueqz}
can be used not only when $\nu=z$, but also in the case where $\nu$ and $z$ are
close, but not exactly equal.  The three solutions given in
Eq.~\eqref{eq:PSD_nueqz} change from paths of steepest ascent to paths of
steepest descent, depending on the sign of $\eta$. In particular, when $\arg
z=0$, then the PSD for $J_\nu(z)$ close to the saddle point $t=0$ satisfies
$\arg t=-2\pi/3$ when $\imag t<0$ and $\arg t=2\pi/3$ when $\imag t>0$.  The
PSD thus approaches the origin from the lower left and departs again to the
upper left [see also Fig.~\ref{fig6}(c)].  

The confluence of the two saddle
points in the case $\nu = z$ represents the major obstacle 
in a numerical treatment of the Bessel function, 
in the context of the saddle-point integration.

%
%

%
%
\subsection{Numerical integration}

In terms of the PSD configurations,
there are three cases that should be considered separately,
depending on the sign of $\imag\left( t_+ \right)$,
where $t_+$ is the saddle point with $\real(t_+) > 0$.
The saddle point configuration
depends on this sign~\cite{SmdHTi2009}.

{\it Case 1:} If $\imag\left( t_+ \right) < 0$ and consequently $\imag\left(
t_- \right) > 0$, the configuration of the saddle points and PSDs in this case
is illustrated in Fig.~\ref{fig6}(a). As is gathered from this plot, one
saddle point contributes to the value of $J_\nu(z)$, and the other to the value
of $H_\nu^{(1)}(z)$. The remaining function $Y_\nu(z)$ can then be computed as
$Y_\nu(z)=-\ii[H_\nu^{(1)}(z)-J_\nu(z)]$.  {\it Case 2:} If $\imag\left( t_-
\right) <0$ and therefore $\imag\left( t_+ \right) > 0$, the situation is as
shown in Fig.~\ref{fig6}(b).  To obtain $J_\nu(z)$, one has to pass both
saddle points, and thus pick up a contribution from both.  {\it Case 3:} The
case $\imag t_+ = \imag t_- =0$ actually needs a separate consideration. It is
shown in Fig.~\ref{fig6}(c), but only in the special case of coalescing
saddle points (see the discussion below). In general, the saddle points may be
separated, still, with $\imag t_+ = \imag t_- =0$. 

For $\nu\approx z$, but not exact equality, the 
two saddle points almost coalesce.
In Fig.~\ref{fig6}(c), we investigate the case $\nu = z$,
with $t_+=t_-=0$. In this case, the saddle point has 
merged together into a double saddle point
at the origin, and $J_\nu(z)$ and
$H^{(2)}_\nu (z)$ are calculated by following the
appropriate paths indicated by red and blue arrows,
respectively. The linear segment approximation
to the true PSD [dashed lines in  Fig.~\ref{fig6}(c)] are seen to be valid
rather far from the saddle point, which implies that the linear approximation
to the PSD may be used to good effect. The given
linear approximation to the PSD for $\nu= z$ may be used also for the case
$\left|1-\frac{\nu}{z} \right| <h$ with given $h$.
In practice, the value of $h$ increases with increasing $|z|$, for a given 
prescribed accuracy of the numerical evaluation. 
As a rule of thumb, the 
imaginary part of $g(t)$ along the effective path should be of order 1, 
where the length of the path is estimated by Eq.~\eqref{eq:epsilon2}. 
 
\begin{figure}[t!]
\includegraphics[width=0.4\textwidth]{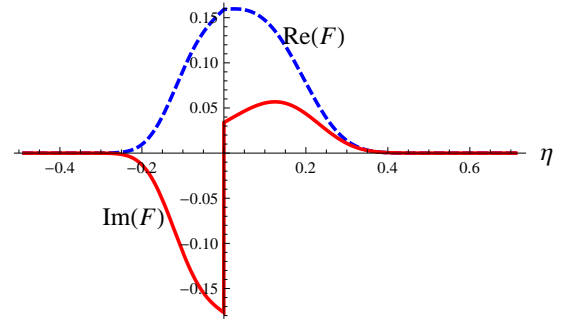}
\caption{\label{fig7} (Color online.)  Real (dashed, blue line) and
imaginary (red, solid line) parts of the integrand $F(\eta)$ from
Eq.~\eqref{eq:integrandF}. The parameter values are 
$\nu=(500 + 1/5)\ee^{3\ii\pi /10}$, and 
$z=(500 + 1/10)\ee^{\ii\pi (3/10+10^{-3})}$. The limits of the
$\eta$-axis correspond to those obtained from the estimate \eqref{eq:epsilon2}}
\end{figure}

To gain a little more insight into the difficulties associated 
with the saddle points approaching each other, we discuss the transition from
Fig.~\ref{fig6}(a) to~\ref{fig6}(c), by a suitable change of the
parameters. The PSD for $J_\nu(z)$ (red arrows) then changes from having smooth
curvature at the saddle point [Fig.~\ref{fig6}(a)] to having a ``kink'' at
the saddle point [Fig.~\ref{fig6}(c)]. Therefore, even if the two saddle
points are separated, at some point the linear approximation to the PSD will
only hold for a very short segment of the true PSD, and is unsuitable. For some
small distance between the saddles, it is instead better to 
pretend that the suitable contour is given by the 
double saddle point and to use the dashed
lines in Fig.~\ref{fig6}(c), even if the saddle points do not exactly
coalesce.  We illustrate the behavior of the integrand along the approximate
PSD when $\nu \approx z$  in Fig.~\ref{fig7}. Shown are real and
imaginary parts of the integrand $F(\eta)$, where
\begin{equation}
J_\nu (z)\approx
 \int_{\eta_{\textrm{min}}}^{\eta_{\textrm{max}}}
F(\eta)\dd \eta,
\end{equation}
\begin{equation}\label{eq:integrandF}
F(\eta)=\frac{1}{2\pi\ii}\ee^{g(t(\eta))} \;
\frac{\dd t(\eta)}{\dd\eta},
\end{equation}
and [see Eq.~\eqref{eq:PSD_nueqz}]
\begin{equation}\label{eq:eta_ex}
t(\eta)=\left\{
\begin{array}{ll}
\eta\left[\cot\left(\dfrac{\pi-\arg z}{3}\right)+\ii\right],&\textrm{if }\eta<0 \,,\\
\eta\left[\cot\left(\dfrac{2\pi-\arg z}{3}\right)+\ii\right],&\textrm{if }\eta\ge 0.
\end{array}\right.
\end{equation}
In Fig.~\ref{fig7}, we have 
$h = \left|1-\frac{\nu}{z} \right| \approx 3\times 10^{-3}$.

%
%

The integration limits for the numerical calculations can be chosen as follows.
The primary objective is to estimate how far out from the saddle point the integration
should be cut off, that is to estimate the limits $\eta_{\textrm{min}}$ and
$\eta_{\textrm{max}}$ in Eq.~ \eqref{eq:approxJnuint}. In the case of separate
saddle points, we expand $g(t)$ around the saddle point up to second order, to
obtain an integral of the form 
\begin{equation}
I_1(X)=\int_0^X \ee^{-C_1 s^2} \dd s \,,
\end{equation}
where $X$ and $C_1$ are positive real constants. It follows that if we want relative
precision $\epsilon=|1-I_1(X)/I_1(\infty)|$ we must cut the integral at
\begin{equation}
\label{eq:epsilon1}
X=\frac{\erfc^{-1}(\epsilon)}{\sqrt{C_1}},
\end{equation}
where $\erfc^{-1}(\cdot)$ is the inverse of the complementary error function.
For coalescing saddle points, also the second derivative of $g(t)$ vanishes,
and we must instead consider integrals like
\begin{equation}
I_2(X)=\int_0^X \ee^{-C_2 s^3} \dd s.
\end{equation}
Here, we should cut at
\begin{equation}\label{eq:epsilon2}
X=\left[\frac{\Gamma^{-1}(\frac{1}{3},\epsilon)}{C_2}\right]^{\frac{1}{3}},
\end{equation}
to obtain a prescribed precision $\epsilon$. In Eq.~\eqref{eq:epsilon2},
$\Gamma^{-1}(\cdot,\cdot)$ is the inverse of the incomplete gamma function with
respect to its second argument.  These estimates assumes a monotonically
decreasing integrand along the path of steepest descent, and it is clear that
it is impossible to reach arbitrary precision due to incurred
oscillations along the linear line segments approximating the 
paths of steepest descent, as one travels too far from the saddle point.

In Ref.~\cite{SmdHTi2009}, the authors also advocate to use an approximation
of the saddle point contour by linear line segments.  They first observe that
it is possible to solve the PSD equation \eqref{eq:PSD} numerically, but then show
that if oscillations do not induce prohibitive numerical instability, equally
accurate results may be obtained by just taking line segments. Because it is
time-consuming to calculate the PSD numerically, they conclude that the
line-segment approach should be favored. However, as mentioned previously, 
this approach necessarily breaks down at some level of precision, especially for large 
values of $|z|$. To reach (in principle) arbitrary precision for large magnitudes 
of the order and argument with the saddle point method, it seems 
unavoidable to follow the true PSD as closely as possible. 

An optimized, hybrid algorithm should perform the
PSD stepping and the integration in parallel: First, one advances on the PSD by one
step, computes the integral on the
segment connecting the previous point by the new point, checks if
the accuracy demand is satisfied, if not, one takes another step. The
step length could be adjusted according to how close the condition
$\imag(g(t_\pm))=\imag(g(t))$ is satisfied. Since
the integrand decreases exponentially along the PSD, {\it and does not oscillate},
it should be straightforward to estimate when to stop the stepping/integration.
Still, for this algorithm to be universally 
applicable, it might be necessary to investigate the 
overlap regions where the saddle-point configurations versus 
the double saddle point configurations should be used,
and to develop special routines that deal with the line 
segments joining the two saddle points in overlapping regions.
We have not pursued this endeavor any further in the current 
work and leave it as an open problem.

\end{document}